\documentclass[]{aimscleaned}

\usepackage{amsmath}
\usepackage{paralist}
\usepackage{graphics} %% add this and next lines if pictures should be in esp format
\usepackage{epsfig} %For pictures: screened artwork should be set up with an 85 or 100 line screen
\usepackage{graphicx}  \usepackage{epstopdf}%This is to transfer .eps figure to .pdf figure; please compile your paper using PDFLeTex or PDFTeXify.
\usepackage[colorlinks=true]{hyperref}
\usepackage{tikz}
\usepackage{lineno}
\usepackage[percent]{overpic}

\hypersetup{urlcolor=blue, citecolor=red}
\usepackage{hyperref}

\def\currentvolume{X}

    \def\ppages{X--XX}
     \def\DOI{10.3934/xx.xx.xx.xx}

\theoremstyle{definition}

\newtheorem{remark}{Remark}

\newcommand{\R}{\mathbb{R}}

\newcommand{\vref}{V_{\textup{ref}}}
\newcommand{\vmax}{V_{\textup{max}}}
\newcommand{\veq}{v^*}
\newcommand{\rhomax}{\rho_{\textup{max}}}
\newcommand{\dmin}{\Delta_{\textup{min}}}
\newcommand{\Dx}{\Delta x}
\newcommand{\Dt}{\Delta t}
\newcommand{\elln}{\ell_{N}}
\newcommand{\DtsuDx}{\lambda}
\newcommand{\nextt}{\textsc{next}}

\newcommand{\Nvcmax}{\Gamma_{\textup{max}}}

\title[An interface-free model for traffic flow] %Use the shortened version of the full title
      {An interface-free multi-scale multi-order model for traffic flow}

\author[Emiliano Cristiani and Elisa Iacomini]{}


 \email{e.cristiani@iac.cnr.it}
 \email{elisa.iacomini@sbai.uniroma1.it}

\thanks{Both authors are members of the INdAM Research group GNCS}

\thanks{$^*$Corresponding author.}

\begin{document}
\maketitle

% Enter the first author's name and address:
\centerline{\scshape Emiliano Cristiani$^*$}
\medskip
{\footnotesize
% please put the address of the first author
 \centerline{Istituto per le Applicazioni del Calcolo ``M. Picone''}
   \centerline{Consiglio Nazionale delle Ricerche}
   \centerline{ Via dei Taurini 19, 00185 Rome, Italy}
} % Do not forget to end the {\footnotesize by the sign }

\vskip 0.5cm%\bigskip

\centerline{\scshape Elisa Iacomini}
\medskip
{\footnotesize
 % please put the address of the second  and third author
 \centerline{ Dipartimento di Scienze di Base e Applicate per l'Ingegneria}
   \centerline{Sapienza, Universit\`a di Roma}
   \centerline{Via Scarpa 16, 00161 Rome, Italy}
}

\bigskip

% The name of the associate editor will be entered by an editorial staff
% "Communicated by the associate editor name" is not needed for special issue.
 %\centerline{(Communicated by the associate editor name)}

%The abstract of your paper
\begin{abstract}
{In this paper we present a new multi-scale method for reproducing traffic flow which couples a first-order macroscopic model with a second-order microscopic model, avoiding any interface or boundary conditions between them. The multi-scale model is characterized by the fact that microscopic and macroscopic descriptions are not spatially separated. On the contrary, the macro-scale is always active while the micro-scale is activated only if needed by the traffic conditions.} The Euler-Godunov scheme associated to the model is conservative and it is able to reproduce typical traffic phenomena like stop \& go waves.
\end{abstract}

\section{Introduction}
In this paper we deal with multi-scale modeling of traffic flow on a single road. 
In particular, we focus on modeling vehicular traffic adopting both a microscopic (agent-based) and a macroscopic (fluid-dynamics) point of view. 
%Such a multi-scale models are nowadays particularly meaningful since both Eulerian (i.e., flux-based) and Lagrangian (i.e., GPS) traffic data are available as real measurements. \rev{A model in which both scales of observation coexist can handle both kind of data without the need of an \emph{ad hoc} preliminary data fusion}. %and can be both used for calibration and validation.
{The main motivation to develop multi-scale models in which both scales of observation coexist is that they can handle both Eulerian (i.e., flux-based) and Lagrangian (i.e., GPS) traffic data, without the need for \emph{ad hoc} preliminary data fusion.}

\subsection{Related work.}
Literature about traffic flow is quite large and a comprehensive review is out of the scope of the paper. For a general introduction we refer the reader to the books by Haberman \cite{habermanbook} and Garavello and Piccoli \cite{piccolibook}, as well as to the survey papers by Helbing \cite{helbing2001RMP} and Piccoli and Tosin \cite{piccoli2012review}.

At a general level, let us just mention the difference between first- and second-order models, see, e.g., \cite{fan2014NHM, fan2013TRR}: The former kind of models represents a simplification of the reality since assume that accelerations are instantaneous, and the traffic conditions are always at equilibrium. 
%In such a models the \emph{velocity} of vehicles is assigned as a function of the (punctual or average) positions of the other vehicles. 
The latter kind of models, instead, is closer to the real dynamics of drivers since considers bounded accelerations. 
%In such a models the \emph{acceleration} of vehicles is assigned as a function of the positions and velocity of the other vehicles. 
Note also that second-order models are able to reproduce some traffic phenomena which are primarily caused by bounded accelerations, like, e.g., stop \& go waves \cite{flynn2009PRE, herty2012KRM, klar2004SIAP, zhao2017TRB}.

\subsubsection*{Many-particle limits.} 
Connections between microscopic and macroscopic traffic flow models are already well established, see, e.g., Ni \cite{ni2011MA} for a short review including meso- and pico-scale. 
Aw et al.\ \cite{aw2002SIAP}, Greenberg \cite{greenberg2001SIAP}, and Di Francesco et al.\ \cite{difrancesco2017MBE} investigated the many-particle limit in the framework of second-order traffic models, deriving the macroscopic Aw-Rascle-Zhang (ARZ) model \cite{aw2000SIAP, zhang2002TRB} from a particular second-order microscopic follow-the-leader (FtL) model \cite{helbing2001RMP, pipes1953JAP}.
Instead Colombo and Rossi \cite{colombo2014RSMUP}, Rossi \cite{rossi2014DCDS-S}, Di Francesco and Rosini \cite{difrancesco2015ARMA}, and Di Francesco et al.\ \cite{difrancesco2017BUMI} investigated the many-particle limit in the framework of first-order traffic models, deriving the macroscopic Lighthill-Whitham-Richards (LWR) model \cite{lighthill1955PRSLA, richards1956OR} as the limit of a first-order FtL model.
Let us also mention the papers by Forcadel et al.\ \cite{forcadel2017DCDS-A}, Forcadel and Salazar \cite{forcadel2015DIE} which investigate the many-particle limit exploiting the link between conservation laws and Hamilton-Jacobi equations.

Moving to road networks, analogous connections are rarer. This is probably due to the fact that macroscopic traffic models on networks are in general ill-posed, since the conservation of the mass is not sufficient on its own to characterize a unique solution at junctions. This ambiguity makes more difficult to find the right limit of the microscopic model, which, in turn, can be defined in
different ways near the junctions. In this context let us mention the paper by Cristiani and Sahu \cite{cristiani2016NHM} which investigates the many-particle limit of a first-order FtL model suitably extended to a road network. The corresponding macroscopic model appears to be the extension of the LWR model on networks introduced by Hilliges and Weidlich \cite{hilliges1995TRB} and then extensively studied by Bretti et al.\ \cite{bretti2014DCDS-S} and Briani and Cristiani \cite{briani2014NHM}.

\subsubsection*{Multi-scale models.} 
The connections between micro- and macro-scale are the foundations of multi-scale models. Such a models couple the traffic description at different scale in order to get a new model which inherits the advantages of the single-scale models.
Many kinds of couplings were proposed so far: first-order FtL and LWR \cite{colombo2015M2AS}, second-order FtL and LWR \cite{garavello2017NoDEA}, second-order FtL and phase-transition model \cite{garavello2013NHM}, second-order FtL and ARZ \cite{lattanzio2010M3AS}. See also, from the engineering literature, the older papers  \cite{bourrel2003TRR, joueiai2015TRR, leclercq2007TRB}. Note that the interface which separates micro- and macro-model can be either fixed or solution-dependent.
In the first case, across the micro-to-macro interface single vehicles (mass particles) are transformed in average density of vehicles, while across the macro-to-micro interface the opposite happens. Obviously, transformation is done preserving the total mass of vehicles.

\subsubsection*{Fundamental diagram}\label{sec:FD}
The fundamental diagram is one of the main ingredients of traffic flow models. It defines the relationship between the flux and the density of vehicles \cite{kerner2002MCM}. It is plain that the flux of vehicles is null in either the case of empty road (null density) or in the case of fully congested road (maximal density, stopped bumper-to-bumper vehicles). For intermediate density levels, real data show a more complicated dynamics. Indeed, drivers act differently in response of the same traffic conditions and, in addition, accelerations and decelerations are far from being instantaneous. 
As a consequence, traffic shows some instabilities \cite{ni2018AMM, wang2013JAT}.

In first-order traffic models, the fundamental diagram can be defined by means of a single function while second-order traffic models allows the fundamental diagram to be multivalued, in the sense that a single value of the density can be associated to many values of the flux, exactly as it happens in reality.
It is interesting to recover the fundamental diagram \emph{a posteriori}, i.e.\ by means of the simulated (predicted) traffic conditions. Ideally, a model should be able to reproduce scattered fundamental diagrams similar to those built on real data. This issue was investigated by many authors, like, e.g., Fan and Seibold \cite{fan2013TRR}, Fan et al.\ \cite{fan2014NHM}, Fan et al.\ \cite{CGARZ_preprint}, Klar et al.\ \cite{klar2004SIAP}, Puppo et al.\ \cite{puppo2016CMS},  Herty and Illner \cite{herty2012KRM}, and  Visconti et al.\ \cite{tosin2017MMS}.

\subsection{Goal.}
This paper proposes a new multi-scale method for reproducing traffic flow which couples a first-order macroscopic model with a second-order microscopic model, avoiding any interface or boundary conditions between them. The macroscopic model \emph{is always running on the whole domain}, while the microscopic model, which can be seen as a ``correction'' of the macroscopic one, appears only when and where the traffic conditions are not in equilibrium. In this way, we make a first-order model like the LWR able to describe second-order effects, and reproduce phenomena like stop \& go waves. At the same time, we avoid to singularly track all the vehicles on the road, which would be computationally expensive.

\subsection{Paper organization}
In Section \ref{sec:singlescale} we recall the classical micro- and macroscopic models which will be used in the rest of the paper.
In Section \ref{sec:cinese} we introduce a new microscopic model derived from the model recently introduced by Zhao and Zhang \cite{zhao2017TRB} to reproduce stop \& go waves. This model will be used as a further ingredient of the multi-scale approach.
In Section \ref{sec:multiscale} we describe the multi-scale model, starting from the basic concepts and then moving to the technical details. We also show that the resulting numerical scheme is conservative (i.e.\ mass-preserving).
In Section \ref{sec:tests} we present some numerical tests.
We conclude the paper with some comments and an overview on future directions.

\section{Classical single-scale models}\label{sec:singlescale}
In this section we recall some classical single-scale models and the relations between them. 

In the microscopic framework, we assume that $N$ vehicles are moving along a single-lane infinite road where overtaking is not possible. 
We denote by $X_k(t)$ the position of the $k$-th vehicle at time $t>0$, for $k=1,\ldots,N$. Similarly, we denote by $V_k(t)$ its instantaneous velocity. 
We assume that all vehicles are equal and that the mass and the length of each of them coincide. We denote by $\elln>0$ the length/mass of each vehicle and by 
$$\mathcal M:=N \elln$$ 
the total mass, which is kept fixed.
We assume that vehicles are ordered in such a way that $X_1(0)<X_2(0)<\ldots<X_N(0)$ and, obviously, $X_{k+1}(0)-X_k(0)\geq\elln$ for all $k=1,\ldots,N-1$. 

In the macroscopic framework, we consider the same mass of vehicles, this time indistinguishable from each other, whose average density at point $x$ and time $t$ will be denoted by $\rho(x,t)\in[0,\rhomax]$ for some maximal density $\rhomax>0$. Conservation of mass will imply that 
$$
\mathcal M=\int_\R \rho(x,t)dx,\qquad t\geq 0.
$$
Similarly, vehicles' average velocity will be denoted by $v(x,t)\in[0,\vmax]$ for some maximal velocity $\vmax>0$.

\medskip

The typical form of a second-order microscopic model of follow-the-leader type is
\begin{equation}\label{FtL2o}
\left\{
\begin{array}{l}
\dot X_k(t)=V_k(t), \quad k\leq N \\ [1mm]
\dot V_k(t)=A(X_k(t),X_{k+1}(t),V_k(t),V_{k+1}(t);\mathbf p),\quad k<N\\ [1mm]
\dot V_N(t)=0
\end{array}
\right.
\end{equation}
where $t>0$, the function $A$ represents the \emph{acceleration} (to be suitably defined) and $\mathbf p$ represents the vector of model parameters. 
Note that the $N$-th vehicle, termed \emph{the leader}, must have a special dynamics since there are no vehicle in front of it. Typically, the leader proceeds at maximal velocity $\vmax$ regardless of the traffic conditions behind. Nonleader vehicles are termed \emph{followers}.

One of the model we will consider in this paper is a well-known and widely used second-order FtL model \cite{aw2002SIAP} defined by
\begin{multline}\label{def:A_ARZmicro}
A\big(X,X',V,V';(\gamma,\tau,\vref,\vmax,\rhomax,\elln)\big)=\\ 
\vref \ \left(\frac{\elln}{\rhomax}\right)^\gamma \ \frac{V'-V}{(X'-X)^{\gamma+1}} + \frac{1}{\tau}\left(\veq\left(\frac{\rhomax\elln}{X'-X}\right)-V\right)
\end{multline}
where
$\gamma\geq 0$, $\tau>0$, $\vref>0$ are additional model parameters, and $\veq:[0,\rhomax]\to[0,\vmax]$ is a $C^1$ decreasing function which gives the equilibrium (i.e.\ desired) velocity of drivers as a function of the degree of congestion. 
For this reason it is natural to assume 
$$\veq(0)=\vmax \quad\text{ and }\quad \veq(\rhomax)=0,$$ 
i.e.\ vehicles move as fast as possible in a free road while completely stop when are bumper-to-bumper. 
%This also implies that $x_{k+1}(t)-x_k(t)\geq\elln$ for all $k$ and $t>0$, provided $\tau$ is sufficiently large \textbf{CONTROLLARE}. 

As recalled in the Introduction, Aw et al.\ \cite{aw2002SIAP} proved that the many-particle limit of the model \eqref{FtL2o}-\eqref{def:A_ARZmicro} is the ARZ model
\begin{equation}\label{ARZ}
\left\{
\begin{array}{ll}
\partial_t\rho+\partial_x(\rho v)=0, & \quad x\in\R,\ t>0 \\ [1mm]
\partial_t(\rho w)+\partial_x(v\rho w)=\rho \frac{\veq(\rho)-v}{\tau}, & \quad x\in\R,\ t>0
\end{array}
\right.
\end{equation}
where
$w(x,t):=v(x,t)+P(\rho(x,t))$ and 
$$
P(\rho):=
\left\{
\begin{array}{ll}
\frac{\vref}{\gamma}\left(\frac{\rho}{\rhomax}\right)^\gamma, & \gamma>0 \\ [2mm]
\vref\ln\left(\frac{\rho}{\rhomax}\right), & \gamma=0.
\end{array}
\right.
$$
Roughly speaking, this means that in the limit $N\to +\infty$ (and $\elln\to 0$) the average density and velocity of particles obeying equation \eqref{FtL2o}-\eqref{def:A_ARZmicro} are given, respectively, by the functions $\rho$ and $v$ solution to equation \eqref{ARZ}.

\medskip

To our purposes, it is important to observe that if vehicles are able to reach instantaneously their equilibrium velocity, 
%how the model \eqref{FtL2o}-\eqref{def:A_ARZmicro} transforms in the case of \emph{traffic equilibrium}. Assuming in \eqref{FtL2o}-\eqref{def:A_ARZmicro} that
%$$
%V_k=V_{k+1}=\veq\left(\frac{\rhomax \elln}{X_{k+1}-X_k}\right), \quad k<N,
%$$
the model \eqref{FtL2o}-\eqref{def:A_ARZmicro} 
drops to first order, simply reading as
\begin{equation}\label{FTL1o}
\dot X_k(t)=\veq\left(\frac{\rhomax \elln}{X_{k+1}-X_k}\right), \quad k<N, \ t>0.
\end{equation}

Similarly to the second-order case, it is possible to compute the many-particle limit of the model \eqref{FTL1o}, which appears to be the classical LWR model
\begin{equation}\label{LWR}
\partial_t\rho+\partial_x(\rho \veq(\rho))=0,\qquad x\in\R,\ t>0,
\end{equation}
(see the Introduction for references).

These relationships suggest and justify the coupling between the models \eqref{FtL2o}-\eqref{def:A_ARZmicro} and \eqref{LWR} which we will detail in Section \ref{sec:multiscale}.

\section{A minimal model for reproducing stop \& go waves.}\label{sec:cinese} 
Stop \& go waves are a typical feature of congested traffic \cite{flynn2009PRE, hoogendoorn2014ITS} and represent a real danger for drivers. Therefore, techniques aimed at reducing them are highly desirable \cite{colombo2004AML,piccoli2018TRC}. 
In this section we describe a microscopic second-order model specifically conceived to reproduce stop \& go waves. This model is nothing but a minimal version of the model recently introduced by Zhao and Zhang \cite{zhao2017TRB} to describe the dynamics of vehicles, bicycles and pedestrians in a unified framework. Our model is 'minimal' in the sense that it is obtained from the Zhao and Zhang's model dropping all the terms which are not strictly necessary to reproduce realistic stop \& go waves. 
It has the form \eqref{FtL2o} with
%\begin{equation}\label{FtL_EE}
%\left\{
%\begin{array}{l}
%\dot X_k(t)=V_k(t),\quad k\leq N \\ [1mm]
%\dot V_k(t)=A(X_k(t),X_{k+1}(t),V_k(t);\mathbf p),\quad k<N\\[1mm]
%\dot V_N(t)=0
%\end{array}
%\right.
%\end{equation}
%where 
\begin{equation}\label{def:Acinese}
A\big(X,X^\prime,V,V^\prime;(\tau,\alpha,\dmin,\vmax)\big)=\frac{1}{\tau}\Big(v^{ZZ}\big(X^\prime-X\big)-V\Big)
\end{equation}
and
\begin{equation}\label{def:v^ZZ}
v^{ZZ}(\Delta):=\left\{
\begin{array}{ll}
0, & \Delta\leq \dmin, \\
\alpha(\Delta-\dmin), & \dmin\leq\Delta\leq\dmin+\vmax/\alpha, \\
\vmax, & \Delta\geq \dmin+\vmax/\alpha.
\end{array}
\right. 
\end{equation}
Here $\alpha>0$ is a parameter and $\dmin>\elln$ is the minimum critical spacing distance between the centers of mass of a vehicle and the preceding one. 
%Therefore, in this case we have $\mathbf p=(\alpha,\dmin,\vmax)$. 
Note that this minimal model, unlike the original one \cite{zhao2017TRB}, is deterministic. Moreover, one should note that the condition $X_{k+1}(0)-X_k(0)\geq\elln \Rightarrow X_{k+1}(t)-X_k(t)\geq\elln$ $\forall t$ is not \textit{a priori} guaranteed.\footnote{The question arises why this condition should hold true in the context of traffic modeling, considering the fact that rear-end collisions are actually possible in real life.}

\begin{figure}[h!]
\centerline{
\includegraphics[width=0.99\textwidth]{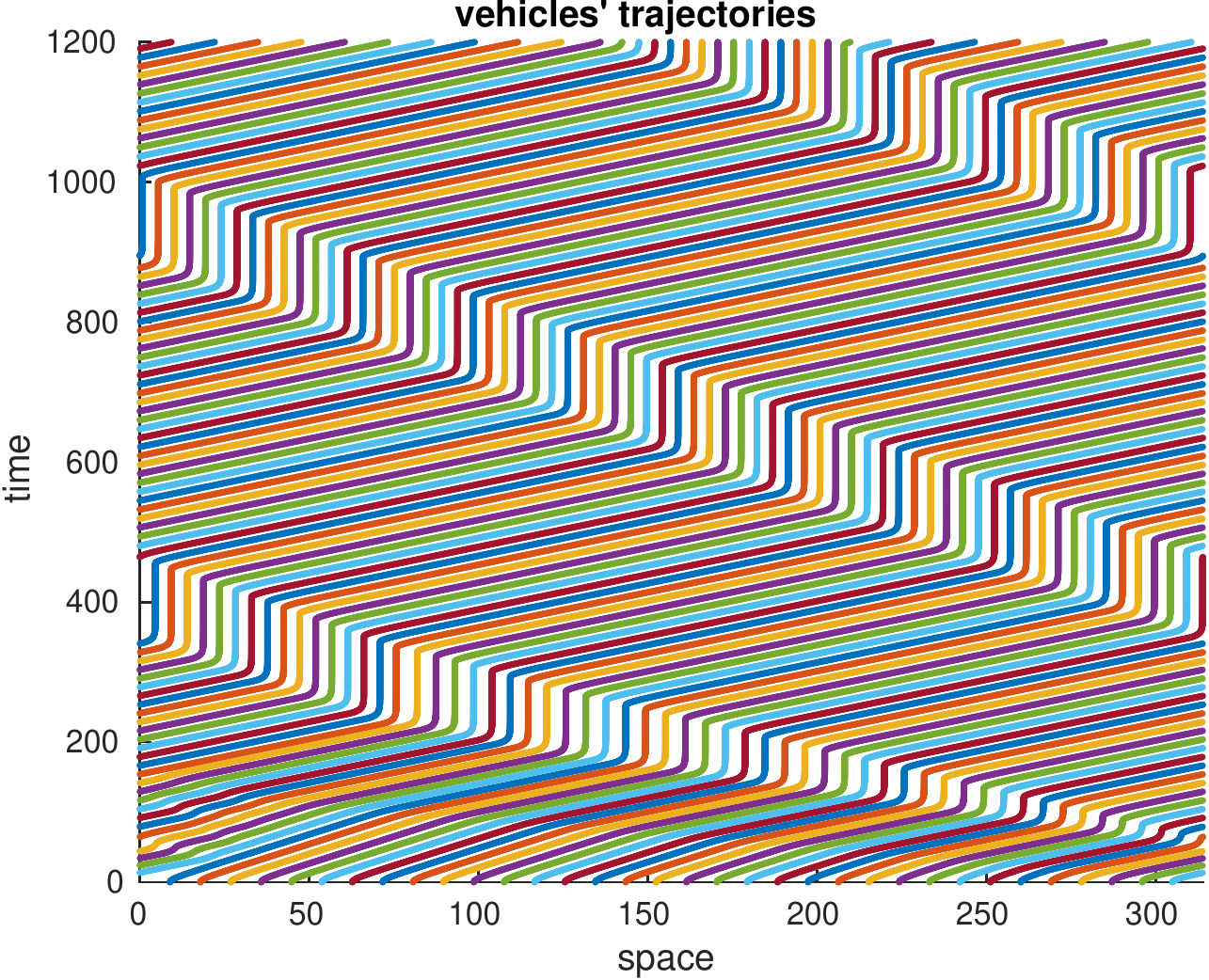}
}
\caption{Space-time trajectories of vehicles obeying to the system \eqref{FtL2o}-\eqref{def:Acinese}-\eqref{def:v^ZZ} with $N=34$, $\alpha=0.6$, $\dmin=7.89$, $\vmax=1$, $\tau=4.86$, $L=314$.}
\label{fig:cinese}
\end{figure}

\begin{figure}[h!]
\centerline{
\includegraphics[width=0.99\textwidth]{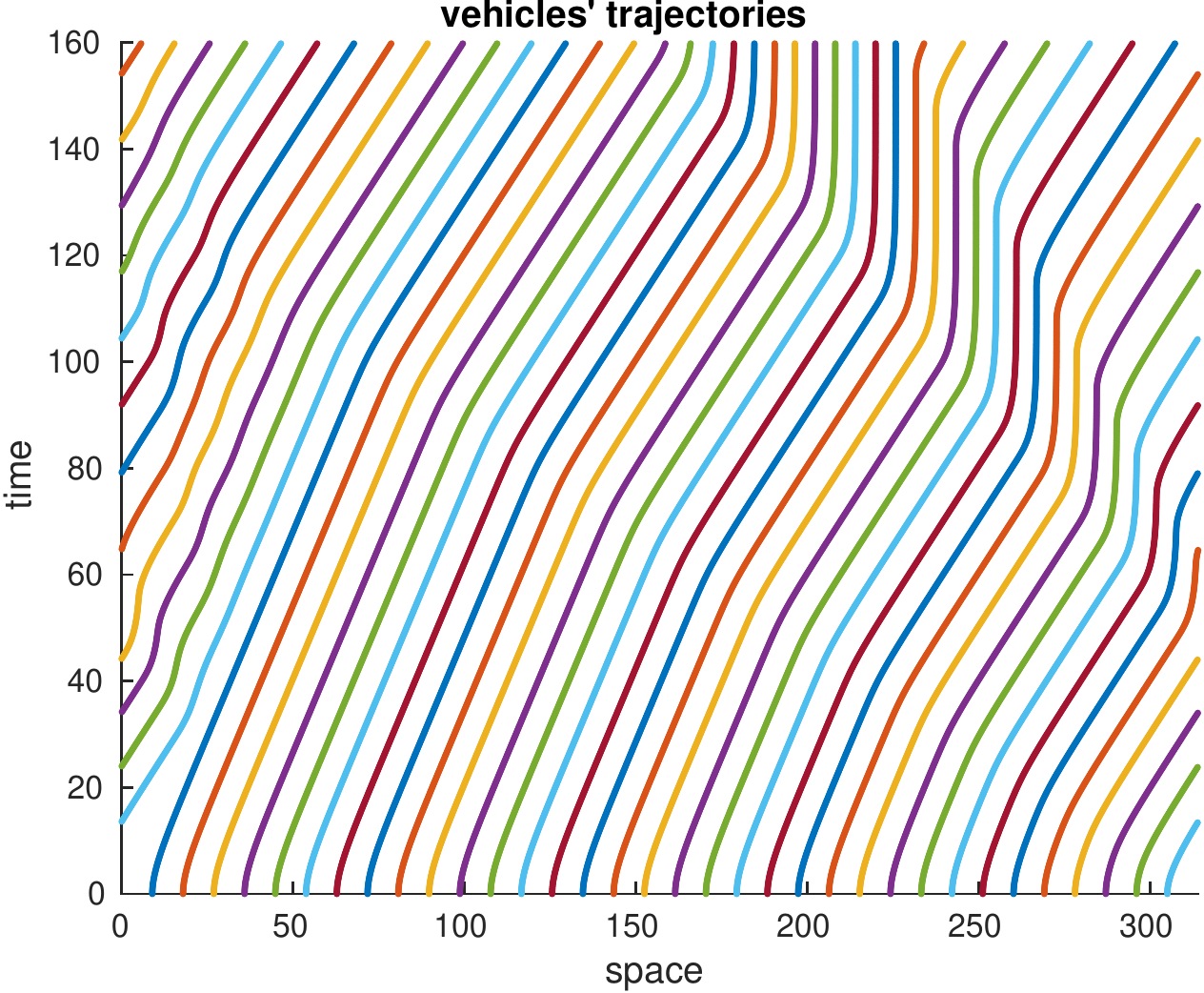}
}
\caption{Zoom of the trajectories shown in Fig.\ \ref{fig:cinese} around initial time. It is well visible the emergence of the stop \& go wave from the interaction between the first and the last vehicle.}
\label{fig:cinese_zoom}
\end{figure}

In Figs.\ \ref{fig:cinese}-\ref{fig:cinese_zoom} we show a typical solution to the system \eqref{FtL2o}-\eqref{def:Acinese}-\eqref{def:v^ZZ} in the case of a \emph{circular road} of length $L$. 
Initial conditions are $X_k(0)=\frac{kL}{N+1}$ and $V_k(0)=0$, for $k=1,\ldots,N$. 
Numerical integration is obtained by the explicit Euler scheme on a road segment $[0,L]$ with periodic boundary conditions. 
It can be seen that backward stop \& go waves are immediately generated by the small perturbation in the initial positions of the vehicles. Indeed, vehicles are initially equispaced
$$
X_k(0)-X_{k-1}(0) =\frac{L}{N+1},\quad k>1
$$ 
with the exception of the couple ($N$,1) (first vehicle in $X_1$ is just in front of the $N$-th vehicle in $X_N$ because of the periodic boundary conditions), for which we have
$$
X_1(0)+L-X_N(0) = \frac{L}{N+1}+L-\frac{NL}{N+1} =\frac{2L}{N+1}.
$$  

Small perturbations in the initial velocity lead to similar effects as well. 

\section{The multi-scale approach}\label{sec:multiscale}
In this section, which is the core of the paper, we describe the multi-scale method, the model and its numerical implementation.

\subsection{General ideas}\label{sec:ideas}
Single-scale models are often unsatisfactory for a number of reasons which involve both modeling and numerical considerations. 
Generally speaking, we can say that:
\begin{itemize}
\item Second-order models are more realistic and often perform better than first-order ones \cite{fan2013TRR, CGARZ_preprint};
\item The numerical approximation of second-order macroscopic models is more difficult than that of first-order macroscopic models, especially if high-order numerical schemes are pursued;
\item Microscopic models require a rather large CPU time if the number of vehicles involved in the simulation is large.
\end{itemize}
Ideally, one would have an easy-to-implement macroscopic model with second-order features. This is the goal which motivated the multi-scale model proposed in this paper.
  
As recalled in the Introduction, multi-scale models are typically based on the \emph{spatial separation} of the microscopic and macroscopic parts. On the contrary, the model we propose here is characterized by the fact that no interface (either fixed or mobile) is explicitly defined. More precisely, the \emph{macroscopic model is always and everywhere alive, while the microscopic model is activated only where and when it is needed}, and it is able to correct (in full or in part) the macroscopic one. 

This procedure is expected to be advantageous if one couples an easy-to-use first-order macroscopic model with a more realistic but still easy-to-use second-order microscopic model. Advantages are complemented by the low computational cost, which comes from the fact that the microscopic model is used only in small parts of the road, therefore the number of vehicles individually tracked is kept low. 

The question arises how the need of the second-order model can be detected. In principle, the second-order model should be activated when and where the traffic is not at equilibrium, i.e.\ the velocity of vehicles is far from the desired one (as for the current traffic conditions). 
This usually happens when nearby vehicles have very different speeds, since this implies the need of  strong accelerations or decelerations which first-order models cannot handle. 
On the contrary, the second-order model can be safely deactivated when vehicles are moving at desired velocity, i.e.\ their acceleration is close to zero.  

%TO DO In addition, macroscopic models cannot take into account the differences among drivers, which are particularly effective in such a situations (different reaction time, etc.).

\subsection{The multi-scale model with complete information}\label{sec:equations}
We are now ready to couple a microscopic model in the form \eqref{FtL2o} 
%(not necessarily with the acceleration \eqref{def:A_ARZmicro}) 
and the macroscopic model \eqref{LWR}. 
In this section we assume, for illustrative purposes, that \emph{both the macroscopic and the microscopic models are alive in the whole space-time domain}, i.e.\ we have a complete information coming from the two models. 

It is given an initial condition at macroscopic level $\rho_0:=\rho(\cdot,0)\in L^1(\R)$, such that $\mathcal M=\int_\R\rho_0(x)dx$ and a number $N$ of microscopic vehicles. 
One can set the initial position of vehicles $\{X_1(0),\ldots,X_N(0)\}$ simply distributing them according to the probability density distribution $\rho_0/\mathcal M$. 
Initial velocities are set as $V_k(0)=\veq(\rho_0(X_k(0)))$, for $k<N$, and $V_N(0)=\vmax$.

Let us also define, as usual, $f:[0,\rhomax]\to\R^+$, $\rho\mapsto f(\rho):=\rho v^*(\rho)$, the flux of vehicles as a function of their density.

The multi-scale model with complete information reads, in integral form, as
\begin{equation}\label{model_analytic}
\left\{
\begin{array}{l}
\displaystyle\partial_t\int_a^b \rho(x,t)dx=
\theta\big(f(a,t)-f(b,t)\big) \\ [1.5mm ]\phantom{xxxx}
\displaystyle
+(1-\theta)\left(
\sum_{k=1}^N\elln\delta(X_k(t)-a)
-\sum_{k=1}^N\elln\delta(X_k(t)-b)
\right),\quad \forall a,b\in\R \\ [2mm]
\dot X_k(t)=V_k(t), \quad k\leq N, \\ [1mm]
\dot V_k(t)=A(X_k(t),X_{k+1}(t),V_k(t),V_{k+1}(t);\mathbf p),\quad k<N,\\ [1mm]
%\dot V_k(t)=\vref \ \elln^\gamma \ \frac{V_{k+1}(t)-V_k(t)}{X_{k+1}(t)-X_k(t)} - \frac{1}{\tau}\left(V_k(t)-\veq\left(\frac{\rhomax\elln}{X_{k+1}(t)-X_k(t)}\right)\right),\quad k<N,\\ [1mm]
\dot V_N(t)=0,
\end{array}
\right.
\end{equation}
where $\theta\in[0,1]$ is an additional parameter, $x\mapsto\delta(x-x_0)$ is the Dirac delta function centered at $x_0$, and the time derivative $\partial_t$ is intended in the distributional sense. 
The idea underlying the model \eqref{model_analytic} is the following: The gain or loss of mass in any space interval $[a,b]$ between time $t$ and $t+dt$ is only given by the flow of vehicles through the boundaries $a$ and $b$, like in classical conservation laws. In this case, at the macroscopic level the flow is given by the classical LWR flux function $f$, while at the microscopic level the mass instantaneously (dis)appears at the passage of vehicles through the boundaries. The parameter $\theta$ is intended for tuning the contribution of the micro- and macro-scale, in the same spirit of \cite{cristiani2015JCSMD, cristiani2011MMS, cristiani2012CDCHawaii, cristianibook, cristiani2018MMS}.
{Here we expect that large values of $\theta$ reduce the oscillations due to abrupt passage of microscopic vehicles across boundaries, while small values of $\theta$ increase the effectiveness of second-order dynamics.}

\subsection{Numerical approximation}\label{sec:numerics}
The system \eqref{model_analytic} can be approximate by a suitable combination of existing numerical schemes. 
In the following, we employ the classical Godunov scheme for the PDE and the explicit Euler scheme for the system of ODEs. 
To do that, we introduce space and time steps $\Dx,\Dt>0$ (with $\DtsuDx:=\frac{\Dt}{\Dx}$) and a grid in space $\{x_j:=j\Dx,\ j\in\mathbb Z\}$ and time $\{t^n:=n\Dt,\ n\in\mathbb N\}$. 
We denote by $C_j:=[x_{j-\frac12},x_{j+\frac12})$ the cell centered in $x_j$. 
As usual, we define 
$$
\rho_j^n:= \frac{1}{\Dx}\int_{C_j}\rho(x,t^n)dx
$$ 
and we denote by $(X_k^n, V_k^n)$ the approximation of $(X_k(t^n), V_k(t^n))$, for $k=1,\dots, N$.

{
In order to define the correspondence between the micro- and the macro-scale we also introduce the scaling parameter $\Nvcmax$, defined as the maximum number of vehicles which can fall in one cell of length $\Dx$. 
$\Nvcmax$ is the microscopic counterpart of $\rhomax$ and it is naturally related to the scaling parameter $\elln$ by 
$$\elln=\frac{\Dx}{\Nvcmax}.$$ %=\frac{L}{N_x\Nvcmax}.$$} 
By means of $\Nvcmax$ we can define the number of vehicles to put in any cell $C_j$, which is equal to $\left\lfloor\frac{\rho_j^0}{\rhomax}\Nvcmax\right\rfloor$. 
In this way we have 
$$
\sum_{k=1}^N \ell_N =
N\ell_N =
\sum_{j\in\mathbb Z} \rho_j^0\Dx =
\int_\R \rho_0(x)dx=\mathcal M.
$$
We also assume that vehicles are initially equispaced in the cell and we assign to each of them the same velocity $v^*(\rho_j^0)$.
This procedure defines the initial positions $\{X_1^0,\ldots,X_N^0\}$ and velocities $\{V_1^0,\ldots,V_N^0\}$ of the microscopic vehicles.
}

The numerical scheme reads as follows.
\begin{equation}\label{scheme}
\left\{
\begin{array}{l}
\displaystyle
\rho_j^{n+1}=\rho_j^n+\theta\DtsuDx\Big(\mathcal G(\rho_{j-1}^n,\rho_{j}^n)-\mathcal G(\rho_j^n,\rho_{j+1}^n)\Big)
+(1-\theta)\DtsuDx
\left(\mathcal F_{j-\frac12}^n-\mathcal F_{j+\frac12}^n\right), \\ [4mm]
X_k^{n+1}=X_k^n+\Dt V_k^n, \quad k\leq N, \\ [3mm]
V_k^{n+1}=V_k^n+\Dt A(X_k^n,X_{k+1}^n,V_k^n,V_{k+1}^n;\mathbf p),\quad k<N,\\ [3mm]
%V_k^{n+1}=V_k^n+\Dt\vref \ \elln^\gamma \ \frac{V_{k+1}^n-V_k^n}{X_{k+1}^n-X_k^n} - \frac{\Dt}{\tau}\left(V_k^n-\veq\left(\frac{\rhomax\elln}{X_{k+1}^n-X_k^n}\right)\right),\quad k<N,\\ [3mm]
V_N^{n+1}=V_N^n,
\end{array}
\right.
\end{equation}
with the classical \textit{Godunov's numerical flux} (see, e.g., \cite{bretti2007ACME, levequebook})
\begin{equation}\label{def:F_G}
\mathcal G(\rho^-,\rho^+):=
\left\{
\begin{array}{ll}
\min\{f(\rho^-),f(\rho^+)\}, & \textrm{if } \rho^-\leq \rho^+ \\
f(\rho^-), & \textrm{if } \rho^->\rho^+ \textrm{ and } \rho^-<\sigma \\
f(\sigma), & \textrm{if } \rho^->\rho^+ \textrm{ and } \rho^- \geq \sigma \geq \rho^+ \\
f(\rho^+), & \textrm{if } \rho^->\rho^+ \textrm{ and } \rho^+>\sigma
\end{array}
\right. 
\end{equation}
(where $\sigma:=\arg\max\limits_{\rho\in[0,\rhomax]}f(\rho)$), and with the \textit{microscopic {numerical} flux}
\begin{equation}\label{def:mathcalF}
\mathcal F_{j\pm\frac12}^n:=\frac{\elln}{\Dt}
\text{Card}\left\{k{\in\{1,\ldots,N\}}:x_k^n<x_{j\pm\frac12}\leq x_k^{n+1}\right\}
\end{equation}
(where 'Card' denote the cardinality of a set).
Note that the microscopic flux is simply computed by counting the number of vehicles passing through the boundaries of the cell $C_j$ in the time interval $[t^n,t^{n+1})$. 

Such a numerical scheme comes along with a natural CFL condition for the microscopic and the macroscopic part, in order to guarantee that in one time step both characteristic curves and vehicles themselves do not move more than one cell apart. To get this, we need to impose
\begin{equation}\label{CFL}
\DtsuDx < \min \left\{ \frac{1}{\max\limits_{\rho\in[0,\rhomax]}|f^\prime(\rho)|},\frac{1}{\vmax} \right\}.
\end{equation}

\subsection{The multi-scale algorithm}\label{sec:algorithm}
One of the main feature of the proposed multi-scale model is that the  microscopic part of the model is activated only where and when the macroscopic model is expected to fail. 
This forces us to introduce two important modifications with respect to the model with complete information.

First, the number $N$ of microscopic (singularly tracked) vehicles can change in time. To deal with that, we denote by $\Gamma^n$ the total number of vehicles at time $t^n$, and by $\Gamma_j^n$ the number of vehicles which, at time $t^n$, fall in the cell $C_j$. 

Second, vehicle $k+1$ is no longer, in general, in front of vehicle $k$. To overcome this issue, we will denote by $\nextt(k)$ the vehicle in front of the $k$-th one, with the convection that leaders have $\nextt=0$.

Finally, we need some additional positive parameters which rule the activation and deactivation of the second-order microscopic model, which will be denoted by $\delta t$, $\delta v$ and $\delta V$. Let us briefly describe their meaning, which will be even clearer after the description of the algorithm.
\begin{itemize}
%\item $\Nvcmax$ is a scaling parameter which gives the maximal number of vehicles which can fall in one cell of length $\Dx$. It is the microscopic counterpart of $\rhomax$ and it is naturally related to $\elln$. 
\item $\delta t$ is the minimal period of time that one microscopic vehicle is active. In other words, if a vehicle is activated at time $t^n$, it cannot be deactivated before time $t^n+\delta t$.
\item $\delta v$ is the minimal variation of the macroscopic velocity function which activates the second-order model.  
\item $\delta V$ is the maximal difference between the current velocity and the equilibrium velocity of microscopic vehicles which allows the deactivation of the second-order model.
\end{itemize}

We are now ready to present the main steps of the algorithm which updates the density and the state of the microscopic vehicles 
$$(\rho^n,X^n,V^n) \to (\rho^{n+1},X^{n+1},V^{n+1}).
$$ 
\begin{enumerate}
\item Compute $\Gamma_j^n$ \ $\forall j$.
\item \emph{Activation of new vehicles (Fig.\ \ref{fig:algoritmo_activation}).} 
For all $j$, if $|v^*(\rho_{j+1}^n)-v^*(\rho_j^n)|>\delta v$, for all $i\in\{j-1, j,j+1, j+2\}$ put new vehicles in cell $C_i$ (unless the cell is already occupied, i.e.\ unless $\Gamma_i^n>0$) with velocity $\veq(\rho_i^n)$. The number of vehicles to put in the cell $C_i$ is $\left\lfloor\frac{\rho_i^n}{\rhomax}\Nvcmax\right\rfloor$. They are initially equispaced in the cell and they have the same velocity.
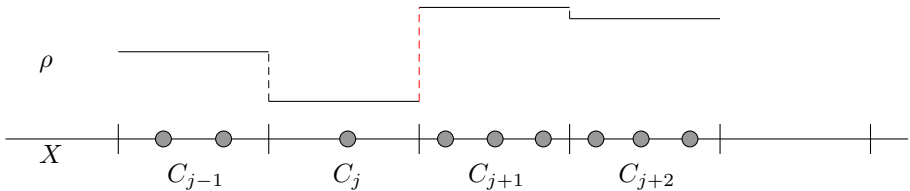
\begin{figure}[h!]
\centering
\begin{tikzpicture}
\draw (1,3) -- (13,3);
\draw (2.5,4.16) -- (4.5,4.16);
\draw (4.5,3.5) -- (6.5,3.5);
\draw (6.5,4.75) -- (8.5,4.75);
\draw (8.5,4.6) -- (10.5,4.6);
\draw[densely dashed,ultra thin] (4.5,3.5) -- (4.5,4.16);
\draw[densely dashed,ultra thin,color=red!100] (6.5,3.5) -- (6.5,4.75);
\draw[densely dashed,ultra thin] (8.5,4.6) -- (8.5,4.75);
\node[text width=11mm] at (2,4) {$\rho$};
\node[text width=11mm] at (2,2.8) {$X$};
\node[text width=11mm] at (3.7,2.5) {$C_{j-1}$};
\node[text width=11mm] at (5.9,2.5) {$C_{j}$};
\node[text width=11mm] at (7.7,2.5) {$C_{j+1}$};
\node[text width=11mm] at (9.7,2.5) {$C_{j+2}$};
%\node[text width=11mm] at (6.75,4.15) {$>\delta v$};
\draw (2.5,2.8) -- (2.5,3.2);
\draw (4.5,2.8) -- (4.5,3.2);
\draw (6.5,2.8) -- (6.5,3.2);
\draw (8.5,2.8) -- (8.5,3.2);
\draw (10.5,2.8) -- (10.5,3.2);
\draw (12.5,2.8) -- (12.5,3.2);
\filldraw[color=black!100,fill=black!40,thin] (3.1,3) circle (3pt);
\filldraw[color=black!100,fill=black!40,thin] (3.9,3) circle (3pt);
\filldraw[color=black!100,fill=black!40,thin] (5.55,3) circle (3pt);
\filldraw[color=black!100,fill=black!40,thin] (6.85,3) circle (3pt);
\filldraw[color=black!100,fill=black!40,thin] (7.51,3) circle (3pt);
\filldraw[color=black!100,fill=black!40,thin] (8.15,3) circle (3pt);
\filldraw[color=black!100,fill=black!40,thin] (8.85,3) circle (3pt);
\filldraw[color=black!100,fill=black!40,thin] (9.45,3) circle (3pt);
\filldraw[color=black!100,fill=black!40,thin] (10.1,3) circle (3pt);
\end{tikzpicture}
\caption{Step 1: Vehicles appear around large jumps of the macroscopic velocity (corresponding to large jumps of the macroscopic density).}
\label{fig:algoritmo_activation}
\end{figure}
\item \emph{Labeling (Fig.\ \ref{fig:algoritmo_labeling}).} 
Find $\nextt(k)$ for all $k$. The rightmost vehicle is labeled as leader ($\nextt=0$). Also, all vehicles $h$ such that $|X_{\nextt(h)}^n-X_h^n|>\Dx$ are also labeled as leader. This choice comes from the assumption that every time a vehicle has a free space of length $\geq\Dx$ in front of it, its dynamics ceases to be dependent on the vehicle in front (if any). We will see, in Step 7, that the dynamics will depend instead on the macroscopic density.\\
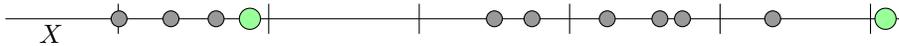
\begin{figure}[h!]
\centering
\begin{tikzpicture}
\draw (1,3) -- (13,3);
\node[text width=11mm] at (2,2.8) {$X$};
\draw (2.5,2.8) -- (2.5,3.2);
\draw (4.5,2.8) -- (4.5,3.2);
\draw (6.5,2.8) -- (6.5,3.2);
\draw (8.5,2.8) -- (8.5,3.2);
\draw (10.5,2.8) -- (10.5,3.2);
\draw (12.5,2.8) -- (12.5,3.2);
\filldraw[color=black!100,fill=black!40,thin] (2.51,3) circle (3pt);
\filldraw[color=black!100,fill=black!40,thin] (3.2,3) circle (3pt);
\filldraw[color=black!100,fill=black!40,thin] (3.8,3) circle (3pt);
\filldraw[color=black!100,fill=green!40,thin] (4.25,3) circle (4pt);
%\filldraw[color=black!100,fill=black!40,thin] (5.1,3) circle (3pt);
%\filldraw[color=black!100,fill=black!40,thin] (5.7,3) circle (3pt);
\filldraw[color=black!100,fill=green!40,thin] (12.7,3) circle (4pt);
\filldraw[color=black!100,fill=black!40,thin] (11.2,3) circle (3pt);
\filldraw[color=black!100,fill=black!40,thin] (10,3) circle (3pt);
\filldraw[color=black!100,fill=black!40,thin] (9.7,3) circle (3pt);
\filldraw[color=black!100,fill=black!40,thin] (9,3) circle (3pt);
\filldraw[color=black!100,fill=black!40,thin] (8,3) circle (3pt);
\filldraw[color=black!100,fill=black!40,thin] (7.5,3) circle (3pt);
\end{tikzpicture}
\caption{Step 3: Green vehicles are leaders.}
\label{fig:algoritmo_labeling}
\end{figure}
\item \emph{Deactivation of followers (Fig.\ \ref{fig:algoritmo_deactivation}).}
Remove all followers $k$ which are active since more than $\delta t$ units of time and such that $\left|V_k-\veq\left(\frac{\rhomax\elln}{X^n_{\nextt(k)}-X_k^n}\right)\right|<\delta V$. Note that, without the first condition new vehicles would immediately deactivated since their velocity is initially at equilibrium (see Step 2). In this way, instead, vehicles have enough time to fully exploit their second-order dynamics. After that, if and when they get close to the equilibrium velocity again, they are deactivated.
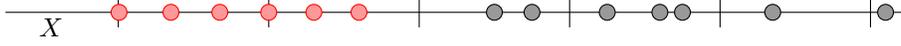
\begin{figure}[h!]
\centering
\begin{tikzpicture}
\draw (1,3) -- (13,3);
\node[text width=11mm] at (2,2.8) {$X$};
\draw (2.5,2.8) -- (2.5,3.2);
\draw (4.5,2.8) -- (4.5,3.2);
\draw (6.5,2.8) -- (6.5,3.2);
\draw (8.5,2.8) -- (8.5,3.2);
\draw (10.5,2.8) -- (10.5,3.2);
\draw (12.5,2.8) -- (12.5,3.2);
\filldraw[color=red!100,fill=red!40,thin] (2.51,3) circle (3pt);
\filldraw[color=red!100,fill=red!40,thin] (3.2,3) circle (3pt);
\filldraw[color=red!100,fill=red!40,thin] (3.85,3) circle (3pt);
\filldraw[color=red!100,fill=red!40,thin] (4.5,3) circle (3pt);
\filldraw[color=red!100,fill=red!40,thin] (5.1,3) circle (3pt);
\filldraw[color=red!100,fill=red!40,thin] (5.7,3) circle (3pt);
\filldraw[color=black!100,fill=black!40,thin] (12.7,3) circle (3pt);
\filldraw[color=black!100,fill=black!40,thin] (11.2,3) circle (3pt);
\filldraw[color=black!100,fill=black!40,thin] (10,3) circle (3pt);
\filldraw[color=black!100,fill=black!40,thin] (9.7,3) circle (3pt);
\filldraw[color=black!100,fill=black!40,thin] (9,3) circle (3pt);
\filldraw[color=black!100,fill=black!40,thin] (8,3) circle (3pt);
\filldraw[color=black!100,fill=black!40,thin] (7.5,3) circle (3pt);
\end{tikzpicture}
\caption{Step 4: Red vehicles are going to be deactivated.}
\label{fig:algoritmo_deactivation}
\end{figure}
\item \emph{Deactivation of leaders.}
Remove all leaders which are not followed by anyone (lonely leaders).
\item Repeat steps 1 and 3 if needed.
\item \emph{Update vehicles' positions and velocities}. We run the microscopic second-order model
\begin{equation}\label{scheme_final_XkVk}
\left\{
\begin{array}{l}
X_k^{n+1}=X_k^n+\Dt V_k^n, \quad \forall k, \\ [2mm]
%V_k^{n+1}=V_k^n+\Dt\vref \ \elln^\gamma \ \frac{V_{\nextt(k)}^n-V_k^n}{X_{\nextt(k)}^n-X_k^n} - \frac{\Dt}{\tau}\left(V_k^n-\veq\left(\frac{\rhomax\elln}{X_{\nextt(k)}^n-X_k^n}\right)\right),\quad \text{if }\nextt(k)>0,\\ [3mm]
V_k^{n+1}=V_k^n+\Dt A(X_k^n, X_{\nextt(k)}^n, V_k^n, V_{\nextt(k)}^n;\mathbf p),\quad \text{if }\nextt(k)>0,\\ [2mm]
V_k^{n+1}=\veq(\rho_{j^*(k,n)+1}^n), \quad\text{if }\nextt(k)=0,
\end{array}
\right.
\end{equation}
where $j^*(k,n)$ is the cell occupied by the vehicle $k$ at time $t^n$. Note that leaders' dynamics only depend on the density $\rho$. Indeed, the velocity of a leader is that of macroscopic vehicles located in the cell in front of it (cf.\ Colombo and Marcellini \cite{colombo2015M2AS}).
\item Compute $\mathcal F^n_{j-\frac12}$ and $\mathcal F^n_{j+\frac12}$ as defined in \eqref{def:mathcalF}.
\item \emph{Update vehicles' density}. We run the multi-scale model, which reads as follows.
For $\theta=0$,
\begin{equation}\label{scheme_final_rho_theta=0}
\rho_j^{n+1}=\rho_j^n+ \DtsuDx \cdot
\left\{
\begin{array}{ll}
\mathcal F^n_{j-\frac12} - \mathcal F^n_{j+\frac12} & \text{ if }\Gamma^n_{j-1},\ \Gamma^n_j,\ \Gamma^n_{j+1}>0 \\ [2mm]
\mathcal F^n_{j-\frac12} - \mathcal G(\rho^n_j,\rho^n_{j+1}) & \text{ if }\Gamma^n_{j-1},\ \Gamma^n_j>0\  \& \ \Gamma^n_{j+1}=0 \\ [2mm]
\mathcal G(\rho^n_{j-1},\rho^n_{j}) - \mathcal F^n_{j+\frac12}& \text{ if }\Gamma^n_{j-1}=0\ \& \ \Gamma^n_{j},\ \Gamma^n_{j+1}>0\\ [2mm]
\mathcal G(\rho^n_{j-1},\rho^n_{j})-\mathcal G(\rho^n_j,\rho^n_{j+1}) & \text{otherwise}.
\end{array}
\right.
\end{equation}
For $\theta\in[0,1]$,
\begin{multline}\label{scheme_final_rho_theta}
\rho_j^{n+1}=\rho_j^n+\\ \DtsuDx \cdot
\left\{
\begin{array}{l}
\left[\theta \mathcal G(\rho_{j-1}^n,\rho_{j}^n)+(1-\theta)\mathcal F^n_{j-\frac12}\right] - \left[\theta \mathcal G(\rho_{j}^n,\rho_{j+1}^n)+(1-\theta)\mathcal F^n_{j+\frac12}\right] \\ 
\phantom{xxxxxxxxxxxxxxxxxxxxxxxxx}
\text{ if }\Gamma^n_{j-1},\ \Gamma^n_j,\ \Gamma^n_{j+1}>0 \\ [2mm]
\left[\theta \mathcal G(\rho_{j-1}^n,\rho_{j}^n)+(1-\theta)\mathcal F^n_{j-\frac12}\right]-\mathcal G(\rho_j^n,\rho_{j+1}^n)  , \\ 
\phantom{xxxxxxxxxxxxxxxxxxxxxxxxx}
\text{ if }\Gamma^n_{j-1},\ \Gamma^n_j>0\  \& \ \Gamma^n_{j+1}=0 \\ [2mm]
\mathcal G(\rho_{j-1}^n,\rho_{j}^n) -\left[\theta \mathcal G(\rho_{j}^n,\rho_{j+1}^n)+(1-\theta)\mathcal F^n_{j+\frac12}\right], \\
\phantom{xxxxxxxxxxxxxxxxxxxxxxxxx}
\text{ if }\Gamma^n_{j-1}=0\ \& \ \Gamma^n_{j},\ \Gamma^n_{j+1}>0\\ [2mm]
\mathcal G(\rho_{j-1}^n,\rho_{j}^n)-\mathcal G(\rho_j^n,\rho_{j+1}^n), \phantom{xxxxx} \text{otherwise}.
\end{array}
\right.
\end{multline}
The scheme is different from the one with complete information \eqref{scheme} because here microscopic information is not always available. Where microscopic vehicles are present, one can choose  a suitable combination of macroscopic and microscopic flux to update the macroscopic density. Where microscopic vehicles are missing, only the macroscopic flux is used, see Fig.\ \ref{fig:algoritmo_updaterho}.
\begin{figure}[h!]
\centering
\begin{tikzpicture}
\shade [fill=blue, fill opacity=0.3] (4.5,3) rectangle (6.5,4.5);
\shade [fill=red, fill opacity=0.4] (6.5,3) rectangle (8.5,3.7);
\shade [fill=green, fill opacity=0.4] (8.5,3) rectangle (10.5,4.2);
\draw (1,3) -- (13,3);
\filldraw[color=black!100,fill=black!40,thin] (4.9,3) circle (3pt);
\filldraw[color=black!100,fill=black!40,thin] (5.5,3) circle (3pt);
\filldraw[color=black!100,fill=black!40,thin] (6.1,3) circle (3pt);
\filldraw[color=black!100,fill=black!40,thin] (6.5,3) circle (3pt);
\filldraw[color=black!100,fill=black!40,thin] (6.8,3) circle (3pt);
\filldraw[color=black!100,fill=black!40,thin] (7.5,3) circle (3pt);
\draw (2.5,2.8) -- (2.5,3.2);
\draw (4.5,2.8) -- (4.5,3.2);
\draw (6.5,2.8) -- (6.5,3.2);
\draw (8.5,2.8) -- (8.5,3.2);
\draw (10.5,2.8) -- (10.5,3.2);
\draw (12.5,2.8) -- (12.5,3.2);
\draw[->] (6.2,4) -- (6.8,4); %freccia
\draw[->] (8.2,4) -- (8.8,4); %freccia
\node[text width=11mm] at (2,4) {$\rho$};
\node[text width=11mm] at (2,2.8) {$X$};
\node[text width=11mm] at (7.9,2.5) {$C_j$};
\node[text width=11mm] at (6.6,4.4) {\footnotesize $\mathcal F_{j-\frac12}$};
\node[text width=11mm] at (8.4,4.4) {\footnotesize $\mathcal G(\rho_j,\rho_{j+1})$};
\end{tikzpicture}
\caption{Step 9: Update of density $\rho_j$ using microscopic flux on the left boundary and macroscopic flux on the right boundary of the cell $j$ (case $\Gamma_{j-1},\ \Gamma_j>0$ \& $\Gamma_{j+1}=0$, $\theta=0$).}
\label{fig:algoritmo_updaterho}
\end{figure}
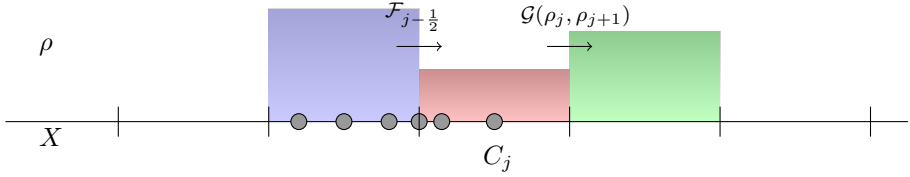
\end{enumerate}

\begin{remark}
We stress once again that the macroscopic density $\rho$ is always updated (Step 9), regardless of the value of $\theta$ and regardless of the presence of microscopic vehicles {(conversely we should manage the interface between the micro- and the macro-scale)}. 
Moreover, the total mass of the system must be evaluated by means on $\rho$ only, as $\int_\R \rho dx$. The appearance of new microscopic vehicles should be seen as a temporary correction procedure which does not imply an additional injection of mass in the system.
Note also that even if $\theta=0$ (dynamics fully driven by the microscopic model), the macroscopic dynamics are still used, precisely where microscopic vehicles are not present. {Moreover the density is used to define the velocity of microscopic leaders}. 
Instead if $\theta=1$ the microscopic model has no effect on the dynamics. \qed 
\end{remark}

\subsection{Conservation of mass}
In the framework of conservation laws and traffic flow models the conservation of the total mass is a crucial property which must be preserved at the discrete level. The multi-scale scheme \eqref{scheme_final_rho_theta} can be easily written in conservative form \cite{levequebook}, therefore it preserves the mass.

In order to reduce the number of cases to be considered, it is convenient to proceed as follows: 
Let us denote by $\mathcal{M}^+_{n,j}$ (resp., $\mathcal{M}^-_{n,j}$) the mass gained (resp., lost) by a generic cell $C_j$ in one time step $n\to n+1$, i.e.,
$$
\Dx\rho^{n+1}_j=\Dx\rho^n_j+\mathcal{M}^+_{n,j}-\mathcal{M}^-_{n,j} \qquad \forall j,\ \forall n,
$$
and then let us simply prove that 
$$
\mathcal{M}^-_{n,j}=\mathcal{M}^+_{n,j+1} \qquad \forall j,\ \forall n.
$$
Let us recall, for the reader convenience, the update rule for $\rho_j^n$ and $\rho^n_{j+1}$ in the case $\theta=0$ (scheme \eqref{scheme_final_rho_theta=0}):
%j
$$
\rho_j^{n+1}=\rho_j^n+ \DtsuDx \cdot
\left\{
\begin{array}{ll}
\mathcal F^n_{j-\frac12} - \mathcal F^n_{j+\frac12} & \text{ if }\Gamma^n_{j-1},\ \Gamma^n_j,\ \Gamma^n_{j+1}>0 \\ [2mm]
\mathcal F^n_{j-\frac12} - \mathcal G(\rho^n_j,\rho^n_{j+1}) & \text{ if }\Gamma^n_{j-1},\ \Gamma^n_j>0\  \& \ \Gamma^n_{j+1}=0 \\ [2mm]
\mathcal G(\rho^n_{j-1},\rho^n_{j}) - \mathcal F^n_{j+\frac12}& \text{ if }\Gamma^n_{j-1}=0\ \& \ \Gamma^n_{j},\ \Gamma^n_{j+1}>0\\ [2mm]
\mathcal G(\rho^n_{j-1},\rho^n_{j})-\mathcal G(\rho^n_j,\rho^n_{j+1}) & \text{otherwise}.
\end{array}
\right.
$$
%j+1
$$
\rho_{j+1}^{n+1}=\rho_{j+1}^n+ \DtsuDx \cdot
\left\{
\begin{array}{ll}
\mathcal F^n_{j+\frac12} - \mathcal F^n_{j+\frac32} & \text{ if }\Gamma^n_{j},\ \Gamma^n_{j+1},\ \Gamma^n_{j+2}>0 \\ [2mm]
\mathcal F^n_{j+\frac12} - \mathcal G(\rho^n_{j+1},\rho^n_{j+2}) & \text{ if }\Gamma^n_{j},\ \Gamma^n_{j+1}>0\  \& \ \Gamma^n_{j+2}=0 \\ [2mm]
\mathcal G(\rho^n_{j},\rho^n_{j+1}) - \mathcal F^n_{j+\frac32}& \text{ if }\Gamma^n_{j}=0\ \& \ \Gamma^n_{j+1},\ \Gamma^n_{j+2}>0\\ [2mm]
\mathcal G(\rho^n_{j},\rho^n_{j+1})-\mathcal G(\rho^n_{j+1},\rho^n_{j+2}) & \text{otherwise}.
\end{array}
\right.
$$
It is easy to see that
$$
\mathcal{M}^-_{n,j}=\mathcal{M}^+_{n,j+1}=\Dt\cdot
\left\{
\begin{array}{ll}
\mathcal F^n_{j+\frac12}, & \text{if }\Gamma^n_{j},\ \Gamma^n_{j+1}>0, \\ [2mm]
\mathcal G(\rho^n_j,\rho^n_{j+1}), & \text{otherwise.}
\end{array}
\right.
$$
Analogously, by \eqref{scheme_final_rho_theta}, we get
$$
\mathcal{M}^-_{n,j}=\mathcal{M}^+_{n,j+1}=\Dt\cdot
\left\{
\begin{array}{ll}
\theta\mathcal G(\rho_j^n,\rho_{j+1}^n)+(1-\theta)\mathcal F^n_{j+\frac12}, & \text{if }\Gamma^n_{j},\ \Gamma^n_{j+1}>0, \\ [2mm]
\mathcal G(\rho^n_j,\rho^n_{j+1}), & \text{otherwise.}
\end{array}
\right.
$$

\section{Numerical tests}\label{sec:tests}
In this section we present some numerical results for the multi-scale model described in Section \ref{sec:algorithm}, with scheme \eqref{scheme_final_rho_theta=0} ($\theta=0$). 
In the first three tests $A$ is defined as in \eqref{def:A_ARZmicro}, while in the last test $A$ is defined as in \eqref{def:Acinese}.

We denote by $L$ the length of the road, by $T$ the final time of the simulation, by $N_x$ the number of space steps, by $N_t$ the number of time steps.
We set $\theta=\gamma=0$, $\vref=\vmax=\rhomax=1$ and $v^*(\rho)=1-\rho$. 
The others parameters used in the simulations are summarized in Table \ref{tab:parameters}. 
\begin{table}[h!]
\begin{tabular}{|c|c|c|c|c|c|c|c|c|c|c|c|}
  \hline 
 & $T$ & $L$ & $N_x$ & $N_t$ & $\tau$ & $\Nvcmax$ & $\delta v$ & $\delta t$ & $\delta V$ & $\alpha$ & $\Delta_{\textup{min}}$
 \\ \hline\hline
T1 & 3     & 20   & 100 & 300   & 0.01      & 20 & 0.08 & 15$\Dt$   & 0.3   & -- & --\\ \hline
T2 & 3     & 20   & 100 & 600   & 0.01--3  & 30 & 0.1   & 15$\Dt$   & 0.5   & -- & --\\ \hline
T3 & 12   & 20   & 100 & 1200 & 0.1        & 30 & 0.1   & 30$\Dt$   & 0.2   & -- & --\\ \hline
T4 & 500 & 314 & 35   & 4000 & 4.86      & 16 & 0.3   & 250$\Dt$ & 0.07 & 0.47 & 2.6$\elln$\\ 
\hline
\end{tabular}
\caption{Model and algorithm parameters used for the numerical tests}
\label{tab:parameters}
\end{table}

\subsection{Test 1 (Activation and deactivation of microscopic model)}
\begin{figure}[h!]
\centerline{
%\textbf{a.}\includegraphics[width=0.47\textwidth]{./fig/T1b_1.pdf} 
%\textbf{b.}\includegraphics[width=0.47\textwidth]{./fig/T1b_100.pdf}
\textbf{a.}
\begin{overpic}[width=0.47\textwidth]{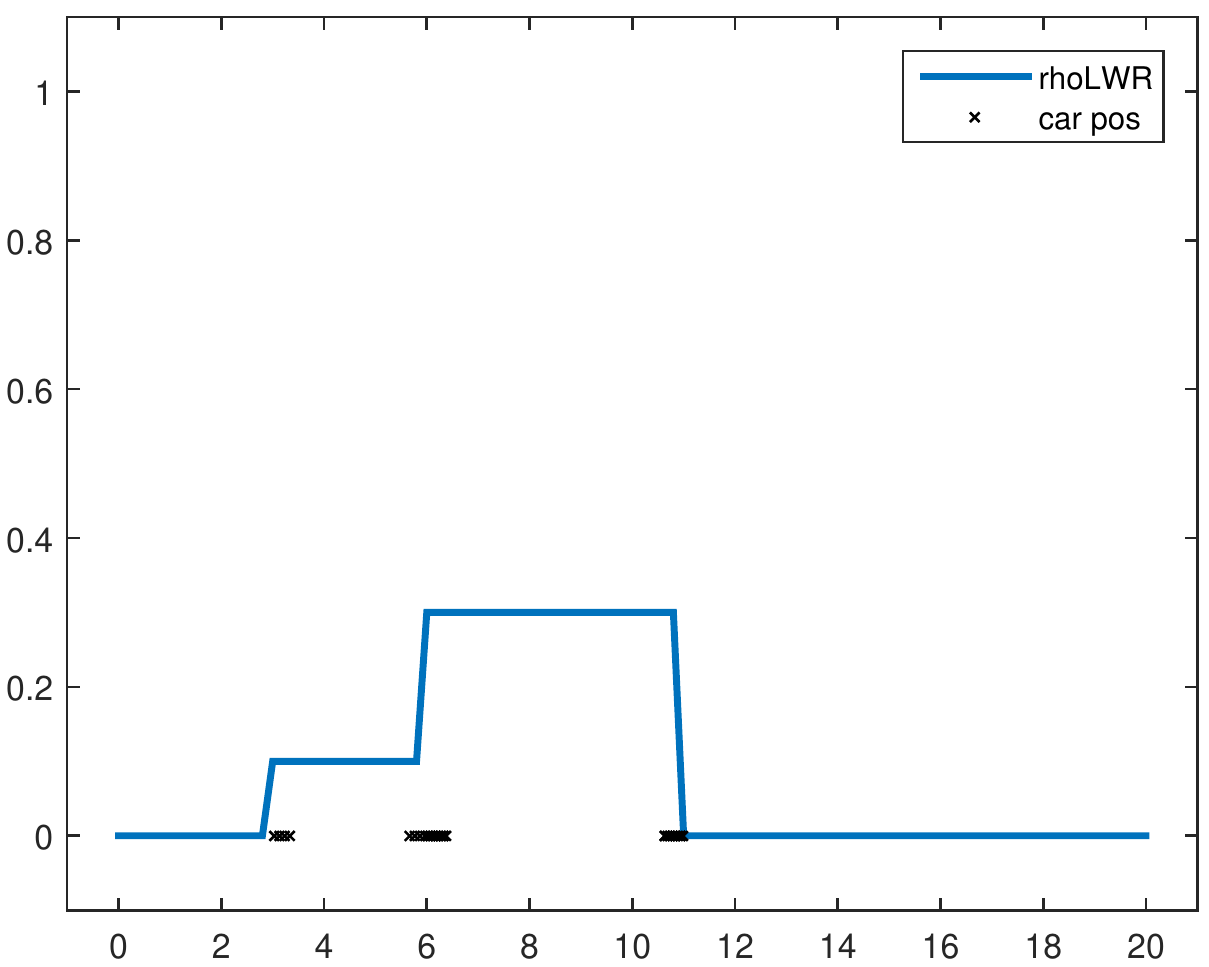}
	\put(52,-4){$x$}\put(-5,45){$\rho$}\end{overpic} 
\textbf{b.}
\begin{overpic}[width=0.47\textwidth]{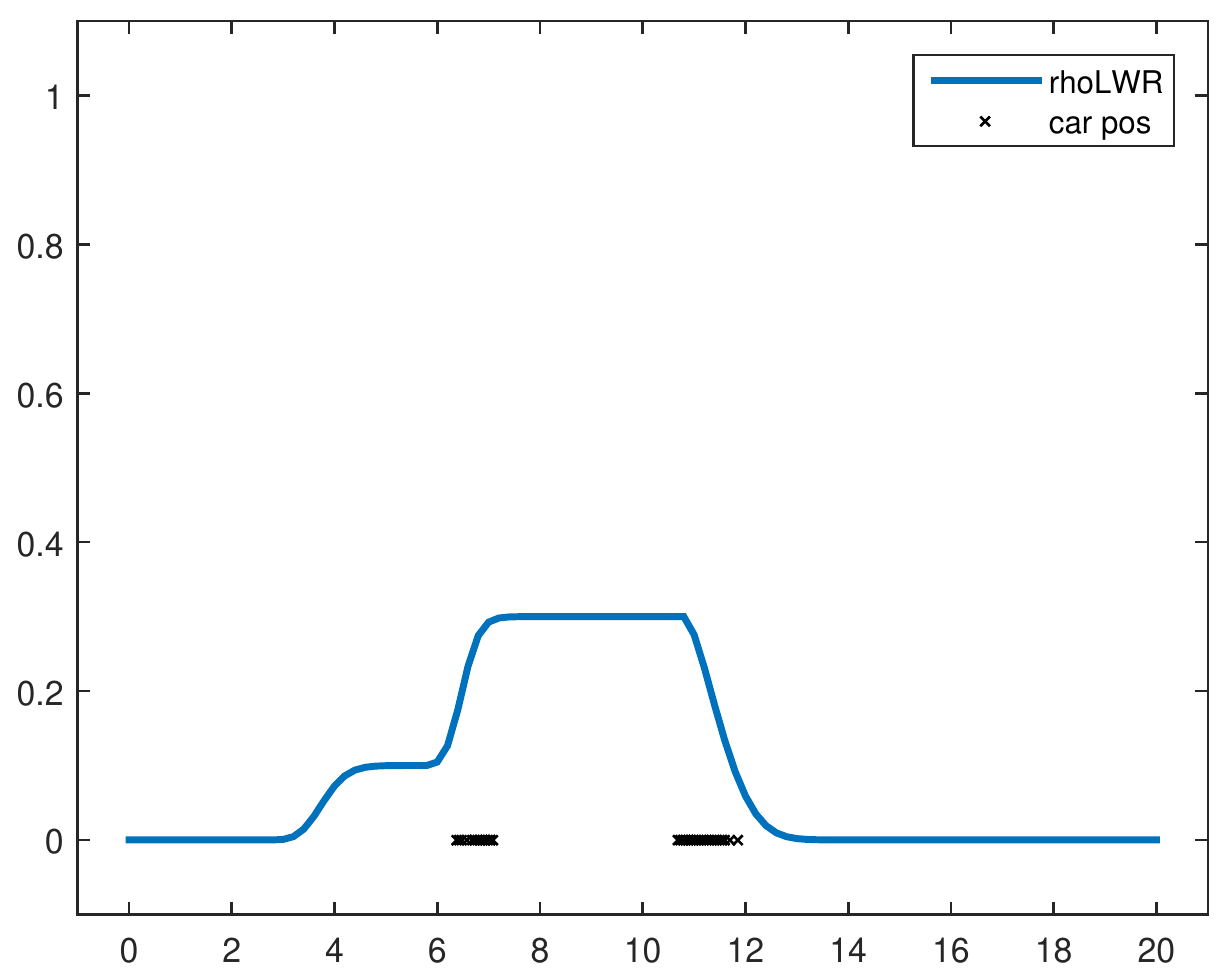}
	\put(52,-4){$x$}\put(-5,45){$\rho$}\end{overpic}  
}
\caption{Test 1: \textbf{a.} $n=1$, \textbf{b.} $n=100$}
\label{fig:T1B}
\end{figure}

In this preliminary test we check the correctness of the steps 2, 3, 4 of the algorithm {and we try to quantify the computational advantage of the new method}. 
To do this, we consider a step function as initial condition $\rho_0$ and we plot the result of the classical LWR model (the density coming from the multi-scale is not plotted for better clarity) and the microscopic vehicles along the $x$-axis. 
We can see that at the initial time vehicles are correctly activated only around the three discontinuities (Fig.\ \ref{fig:T1B}a). After some time, the first discontinuity is smoothed enough for allowing the deactivation of vehicles (Fig.\ \ref{fig:T1B}b).

\begin{figure}[h!]
	\centerline{
		\includegraphics[width=0.6\textwidth]{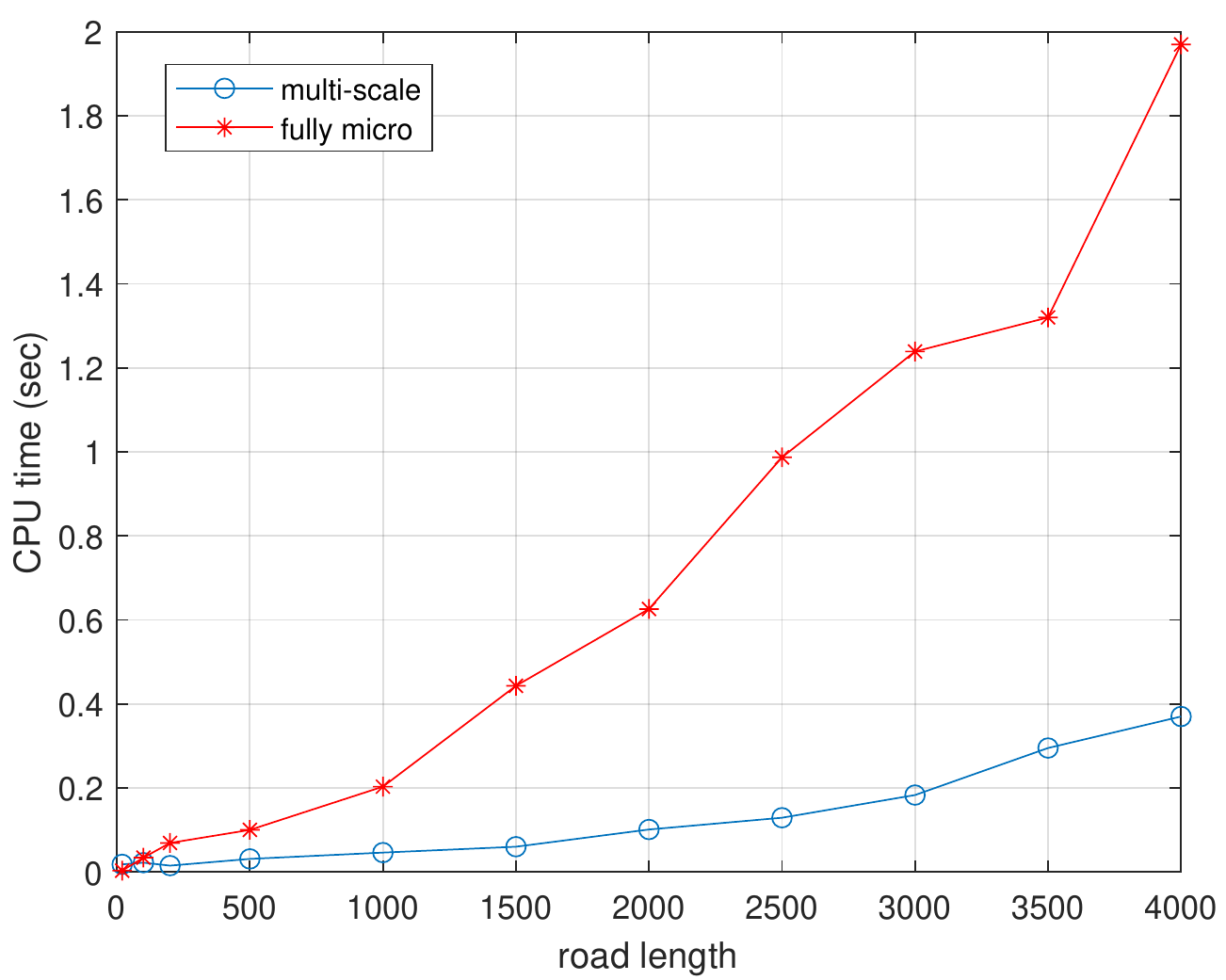}
	}
	\caption{Test 1: CPU time for the fully microscopic model and the multi-scale model.}
	\label{fig:cpu_time}
\end{figure}

{In order to quantify the computational advantage, we compare the CPU time consumed by the fully microscopic model and the multi-scale model (the fully macroscopic model is not considered since we are assuming that it is not able alone to provide a satisfactory traffic description). 
We run the two models for different values of the road length $L$, keeping fixed both $\Dx$ and $\Nvcmax$. Positions of the discontinuities of $\rho_0$ (at $x=3,6,11$ if $L=20$) are scaled linearly with $L$. 
Note that the increase of $L$ (and then of $N_x$ accordingly), causes the increase of the total number of vehicles in the fully microscopic model, while in the multi-scale model the number of vehicles remains constant, being vehicles confined around discontinuities.
In Fig.\ \ref{fig:cpu_time} we show the CPU time for the two models as a function of $L$.}
{We see that the CPU time scales almost linearly with respect to $L$ in both cases. This is a bit surprising for the multi-scale model but it can be explained by the fact that the computational effort of Steps 1, 2, 9 of the algorithm scales linearly with the number of cells $N_x$. 
Apart from that, we see that the CPU time for the multi-scale model is much lower than that of the fully microscopic model, and this is mainly due to the fact that the number of microscopic vehicles is kept low.}

\subsection{Test 2 (effect of the relaxation term $\tau$)}
\begin{figure}[h!]
\centering
\textbf{a.}
\begin{overpic}[width=0.45\textwidth]{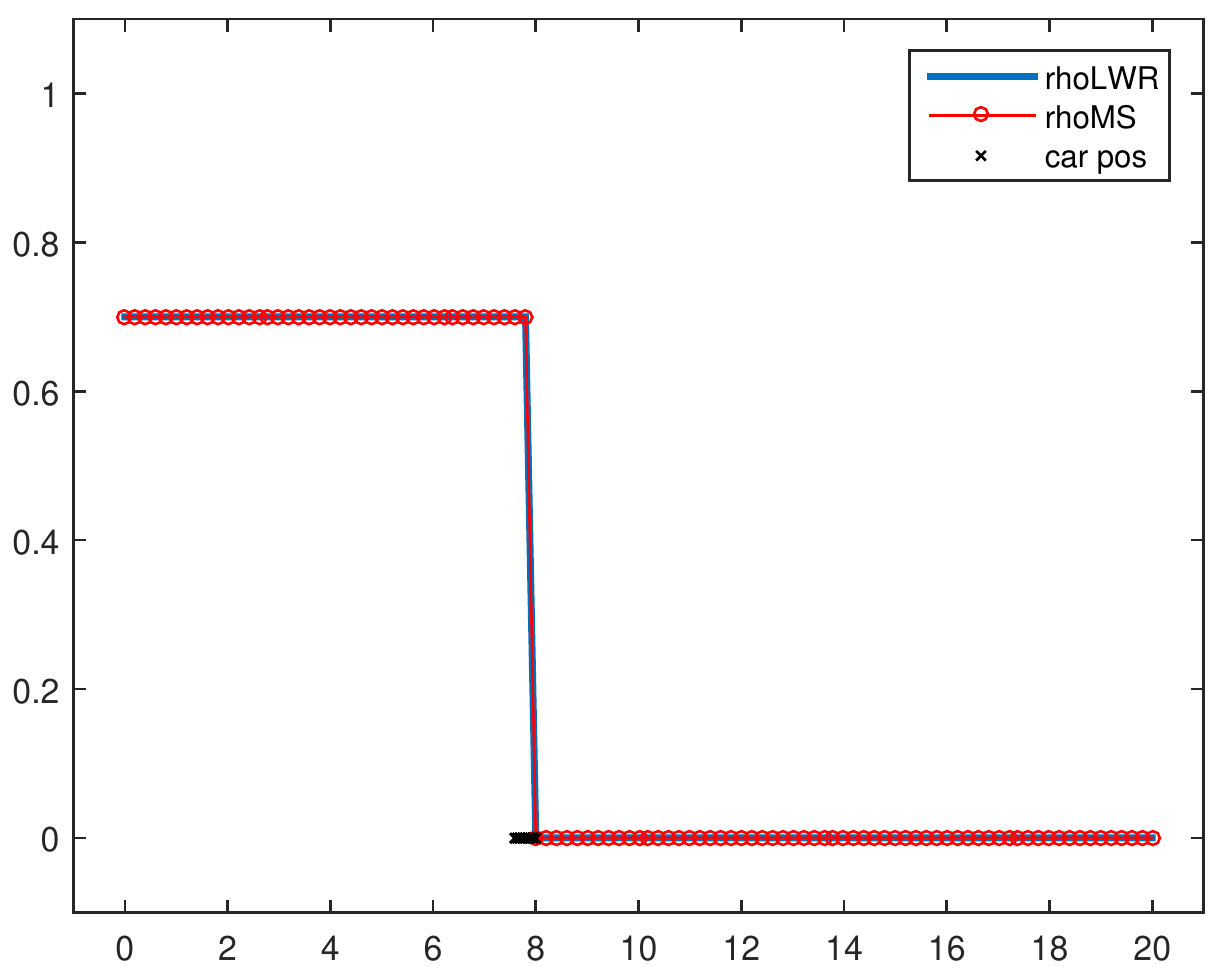}
	\put(52,-4){$x$}\put(-5,45){$\rho$}\end{overpic} 
\textbf{b.}
\begin{overpic}[width=0.45\textwidth]{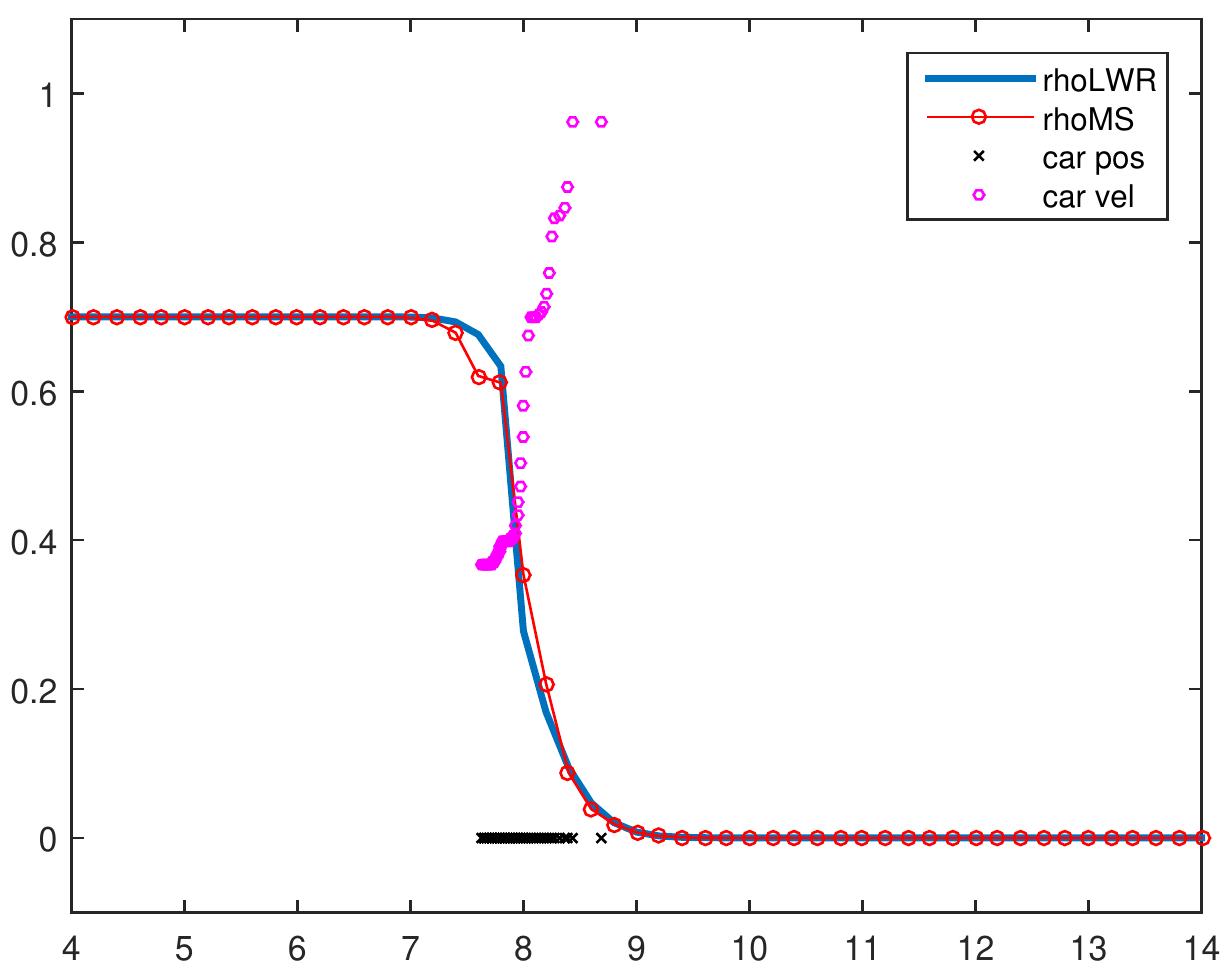}
	\put(52,-4){$x$}\put(-5,45){$\rho$}\end{overpic}
\vskip0.5cm
\textbf{c.}
\begin{overpic}[width=0.45\textwidth]{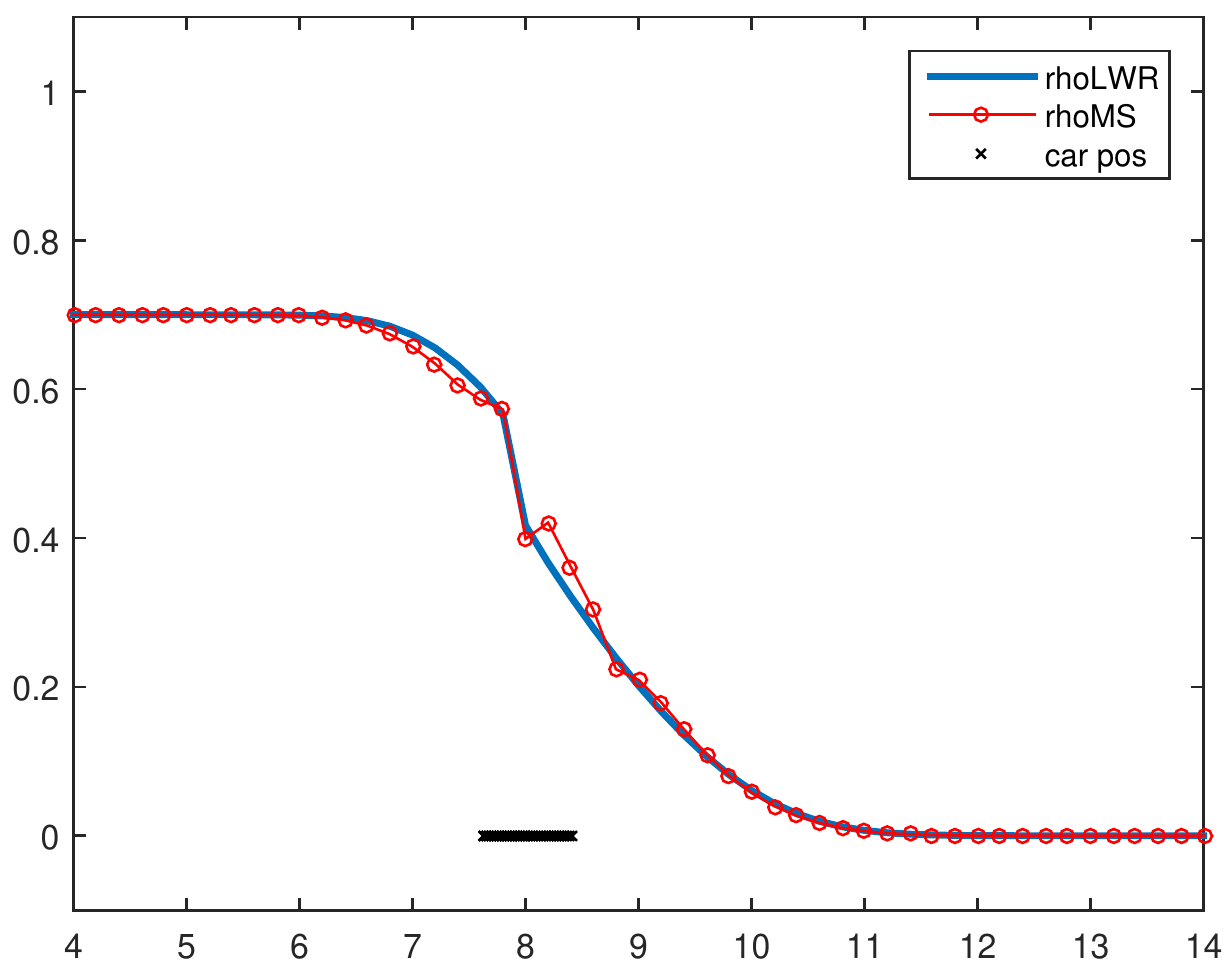}
	\put(52,-4){$x$}\put(-5,45){$\rho$}\end{overpic}
\textbf{d.}
\begin{overpic}[width=0.45\textwidth]{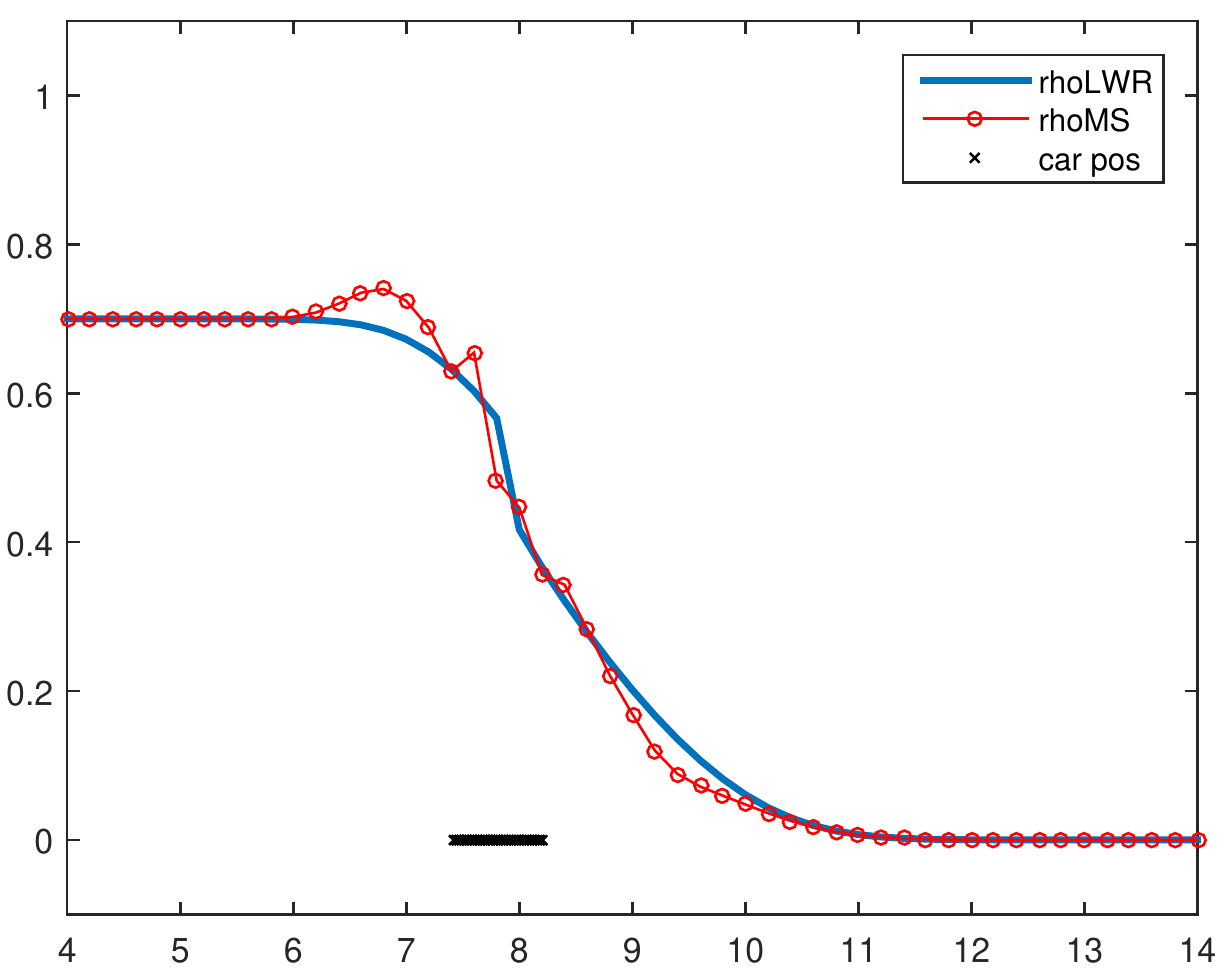}
	\put(52,-4){$x$}\put(-5,45){$\rho$}\end{overpic}
\caption{Test 2: \textbf{a.} $n=1$, \textbf{b.} $n=100$, \textbf{c.} $n=400$, $\tau=0.01$, \textbf{d.} $n=400$, $\tau=3$.}
\label{fig:T2}
\end{figure}
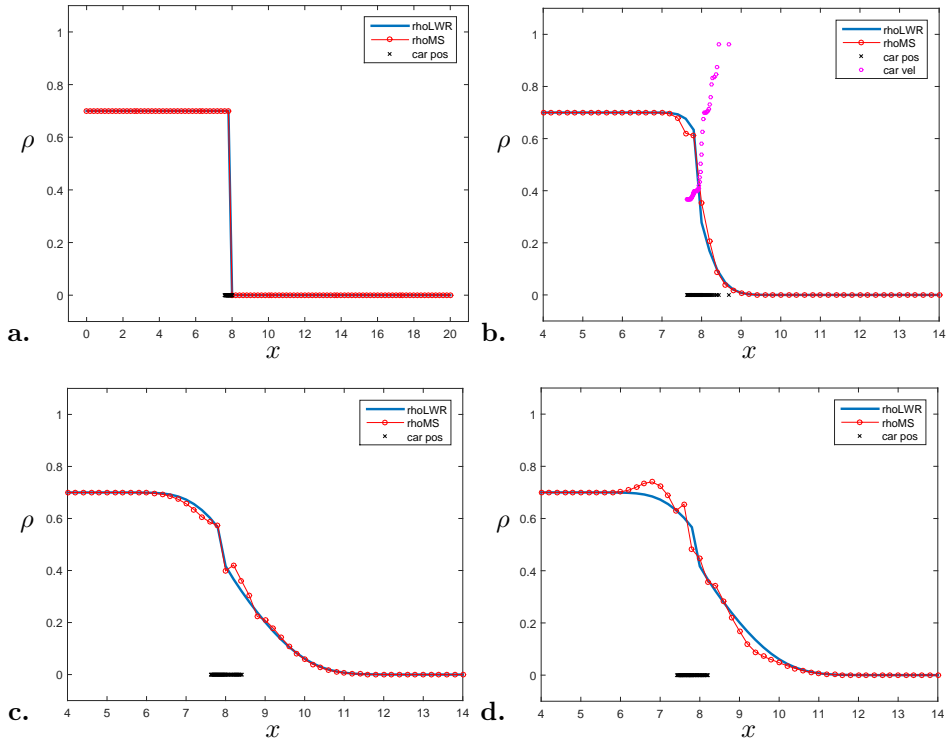

In this test we investigate the effect of the relaxation parameter $\tau$ in the second-order model \eqref{def:A_ARZmicro}. This parameter is related to the \emph{reactivity} of drivers. More precisely, the smaller $\tau$, the more reactive the drivers are and the more the vehicles are able to accelerate and reach the equilibrium velocity rapidly. 
We consider the case of a road congested in the first part and totally free in the second part (Fig.\ \ref{fig:T2}a). 
In this case the LWR model shows immediately the classical rarefaction fan around the original discontinuity (Fig.\ \ref{fig:T2}bcd). 
Microscopic vehicles are activated only around the discontinuity (Fig.\ \ref{fig:T2}abcd) and are able to take into account the bounded acceleration of vehicles, see especially Fig.\ \ref{fig:T2}b where the velocity of microscopic vehicles is also plotted. 
At time step $n=400$ the difference between $\tau=0.01$ (Fig.\ \ref{fig:T2}c) and $\tau=3$ (Fig.\ \ref{fig:T2}d) is quite visible. In the former case (highly reactive drivers) the dynamics of the multi-scale model are very similar to that of the LWR model, while in the latest case (poorly reactive drivers) the multi-scale model differs from LWR and correctly take into account a delay in moving forward. 
%This causes a little queue formed upstream. 
%Note that the parts of the roads with constant density do not show second-order effects and microscopic vehicles are not needed.
%

\begin{figure}[h!]
\centering
\textbf{a.}
\begin{overpic}[width=0.45\textwidth]{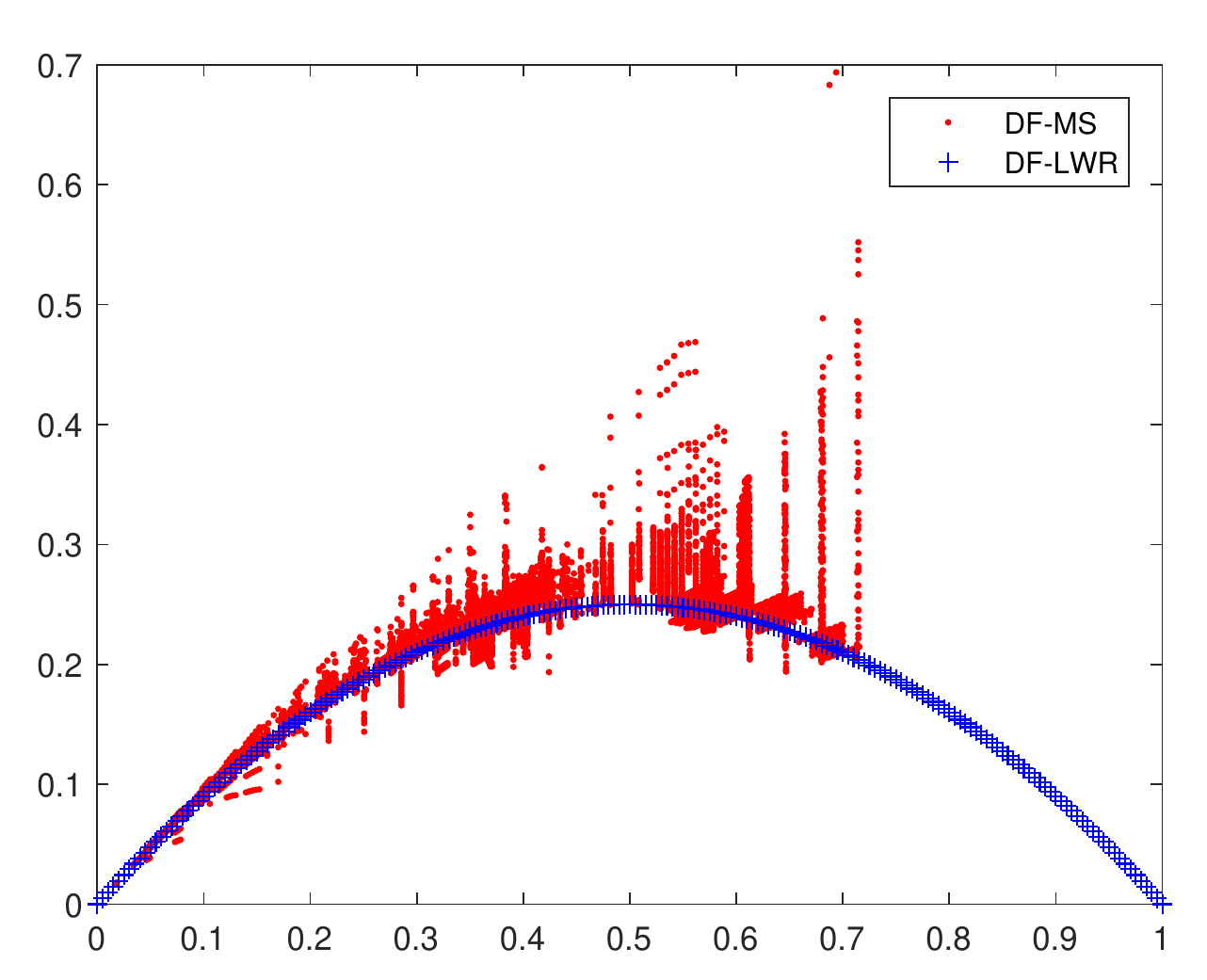}
	\put(52,-3){$\rho$}\put(-6,45){$\rho V$}\end{overpic}
\textbf{b.}
\begin{overpic}[width=0.45\textwidth]{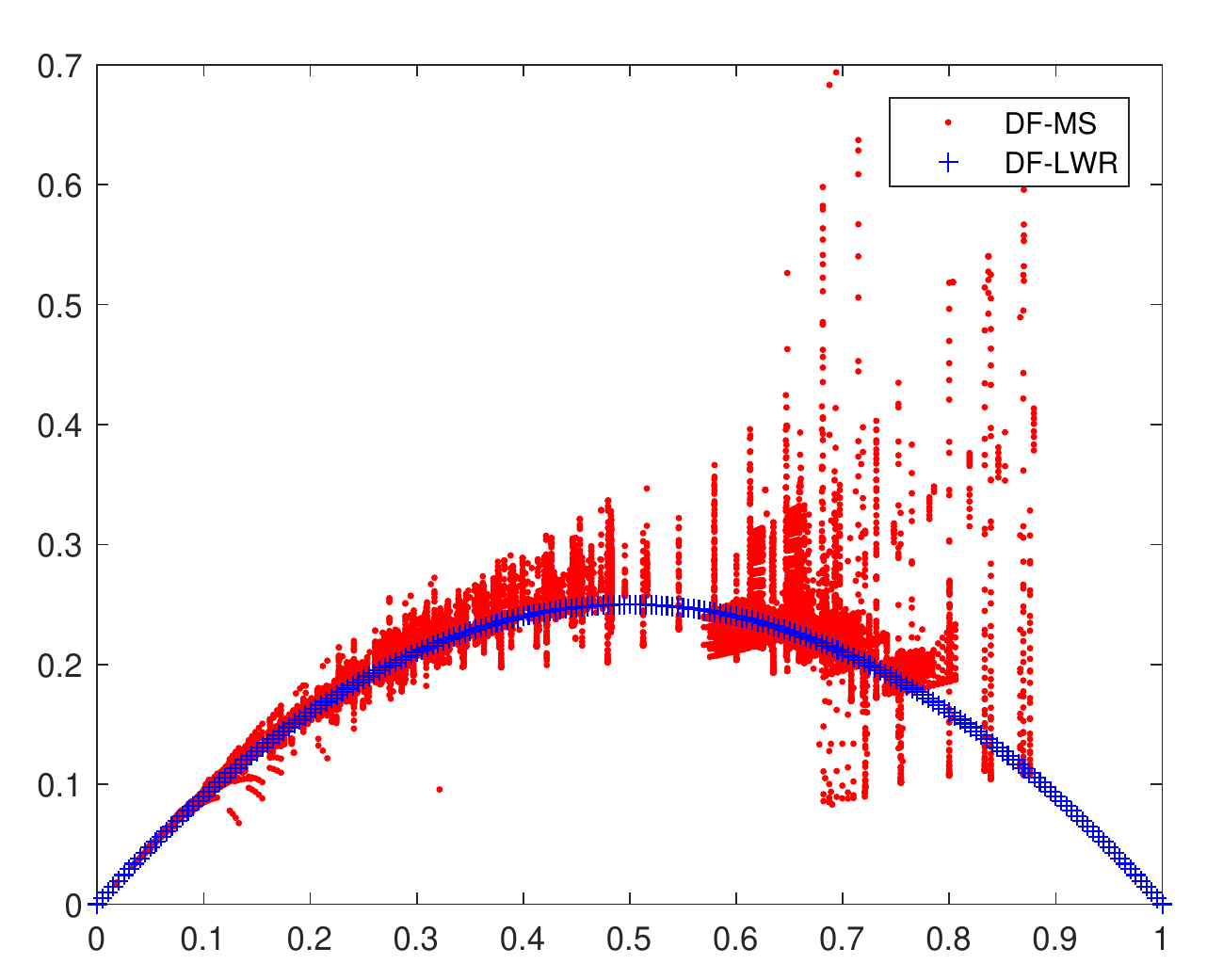}
	\put(52,-3){$\rho$}\put(-6,45){$\rho V$}\end{overpic}
\caption{Test 2: Fundamental diagram of the multi-scale model compared with that of the LWR model. \textbf{a.} $\tau=0.01$, \textbf{b.} $\tau=3$.}
\label{fig:dftest2}
\end{figure}

As discussed in Section \ref{sec:FD}, it is interesting to recover the fundamental diagram \emph{a posteriori}, i.e.\ by means of the simulated traffic conditions. To do that, we plot the set of 2D points
\begin{equation}\label{FDformula}
\{(\rho_{j^*(k,n)}^n, \rho_{j^*(k,n)}^n V_k^n) \ \text{$\forall k$ and $\forall n$, provided vehicle $k$ is active at $t=t^n$}\}
\end{equation}
where $j^*(k,n)$ is the cell occupied by the vehicle $k$ at time $t^n$, see Fig.\ \ref{fig:dftest2}.
It can be seen that the multi-scale model is indeed able to recover a scattered (multivalued) fundamental diagram. Moreover, the scattering increases with $\tau$, as expected, {and decreases smoothly as $\rho\to 0$}.

\subsection{Test 3 (self-sustained perturbation)}
In this test we show the behavior of the model in presence of a perturbation. The perturbation is represented, at time $t=0$, by a small region where vehicles are moving slower than elsewhere (and therefore their density is higher). A typical example is given by \emph{sags}, which are road sections along which gradient changes significantly from downwards to upwards \cite{hoogendoorn2014ITS}. 

\begin{figure}[h!]
	\centering
	\textbf{a.}
	\begin{overpic}[width=0.45\textwidth]{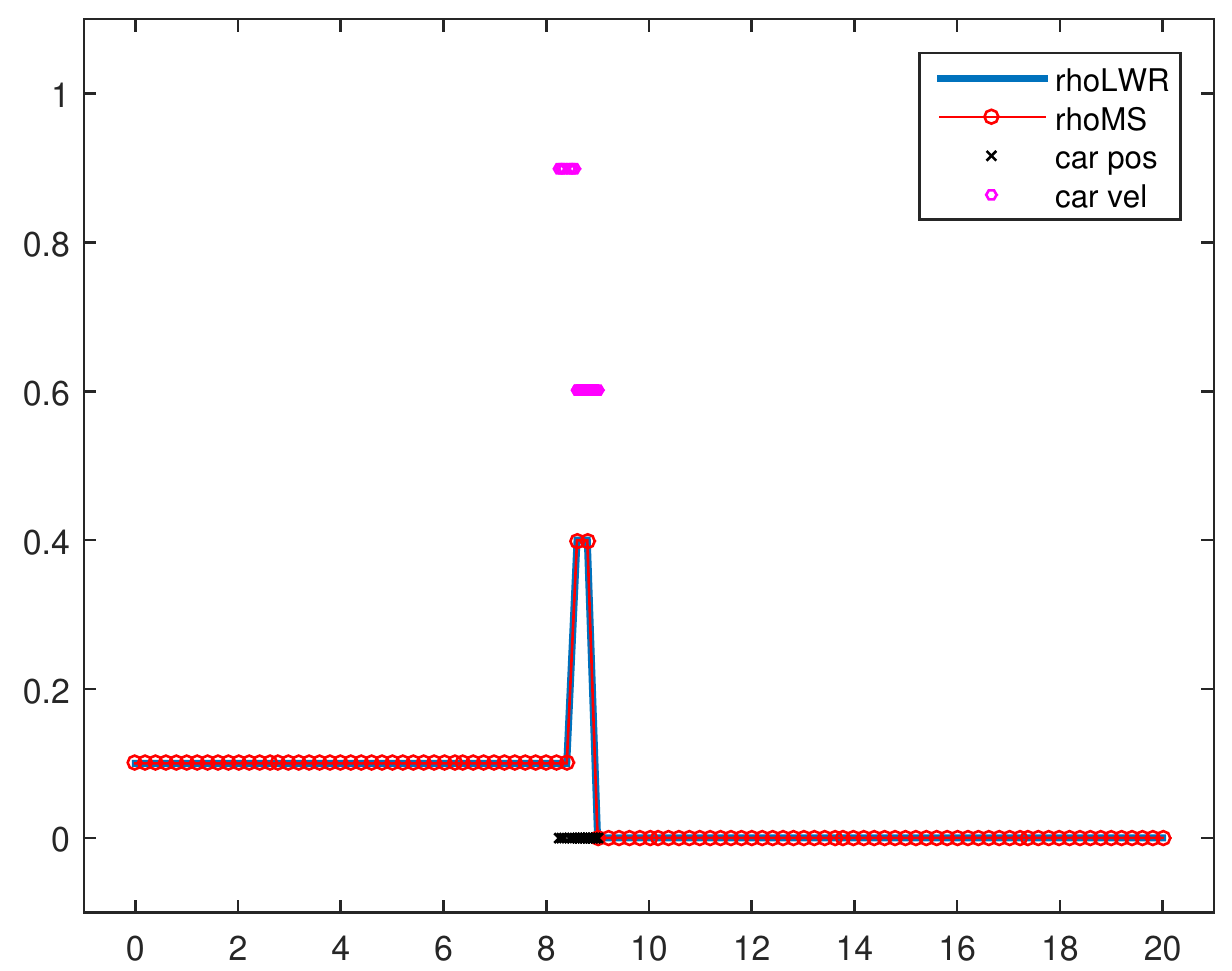}
		\put(52,-4){$x$}\put(-5,45){$\rho$}\end{overpic} 
	\textbf{b.}
	\begin{overpic}[width=0.45\textwidth]{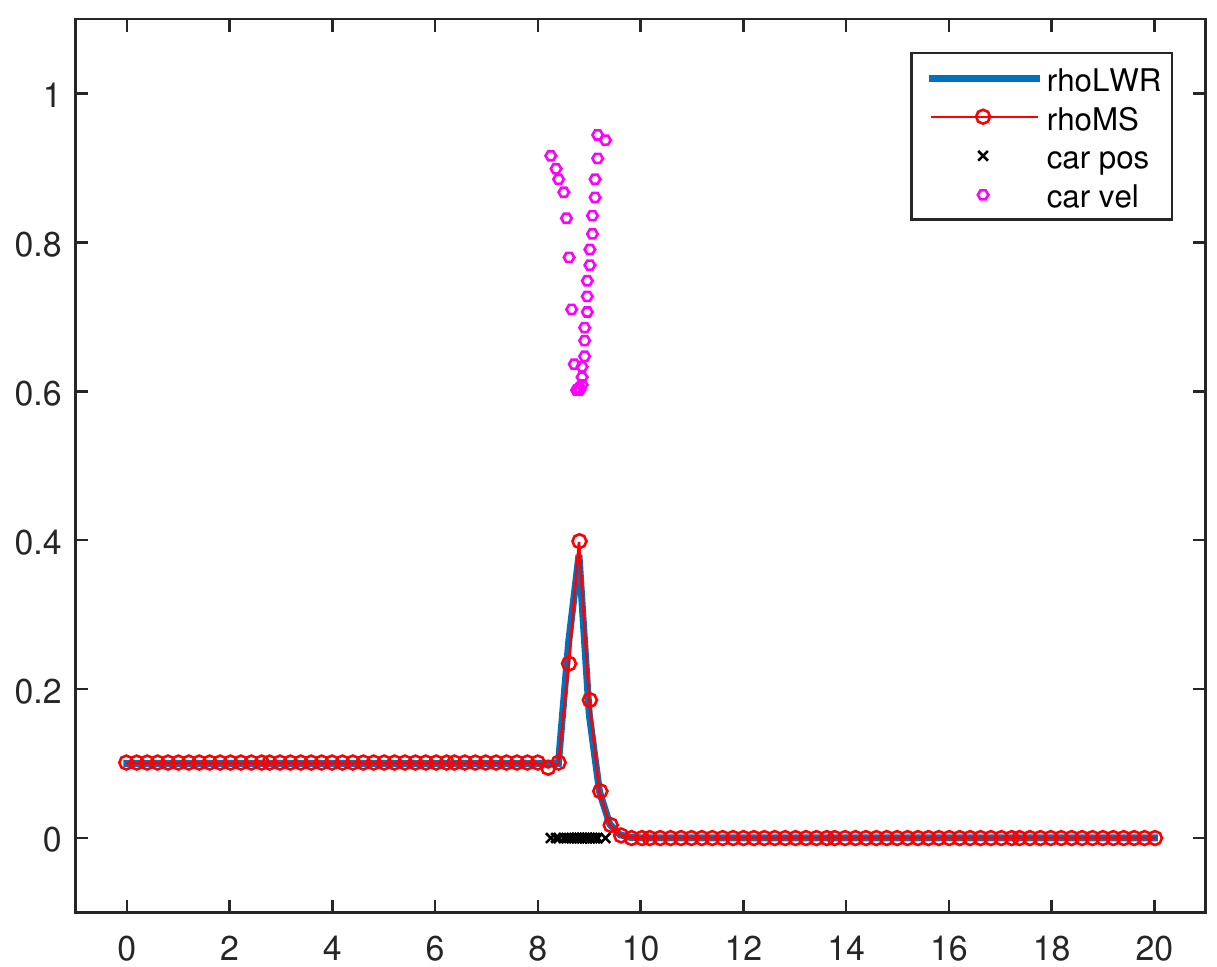}
		\put(52,-4){$x$}\put(-5,45){$\rho$}\end{overpic}
	\vskip0.5cm
	\textbf{c.}
	\begin{overpic}[width=0.45\textwidth]{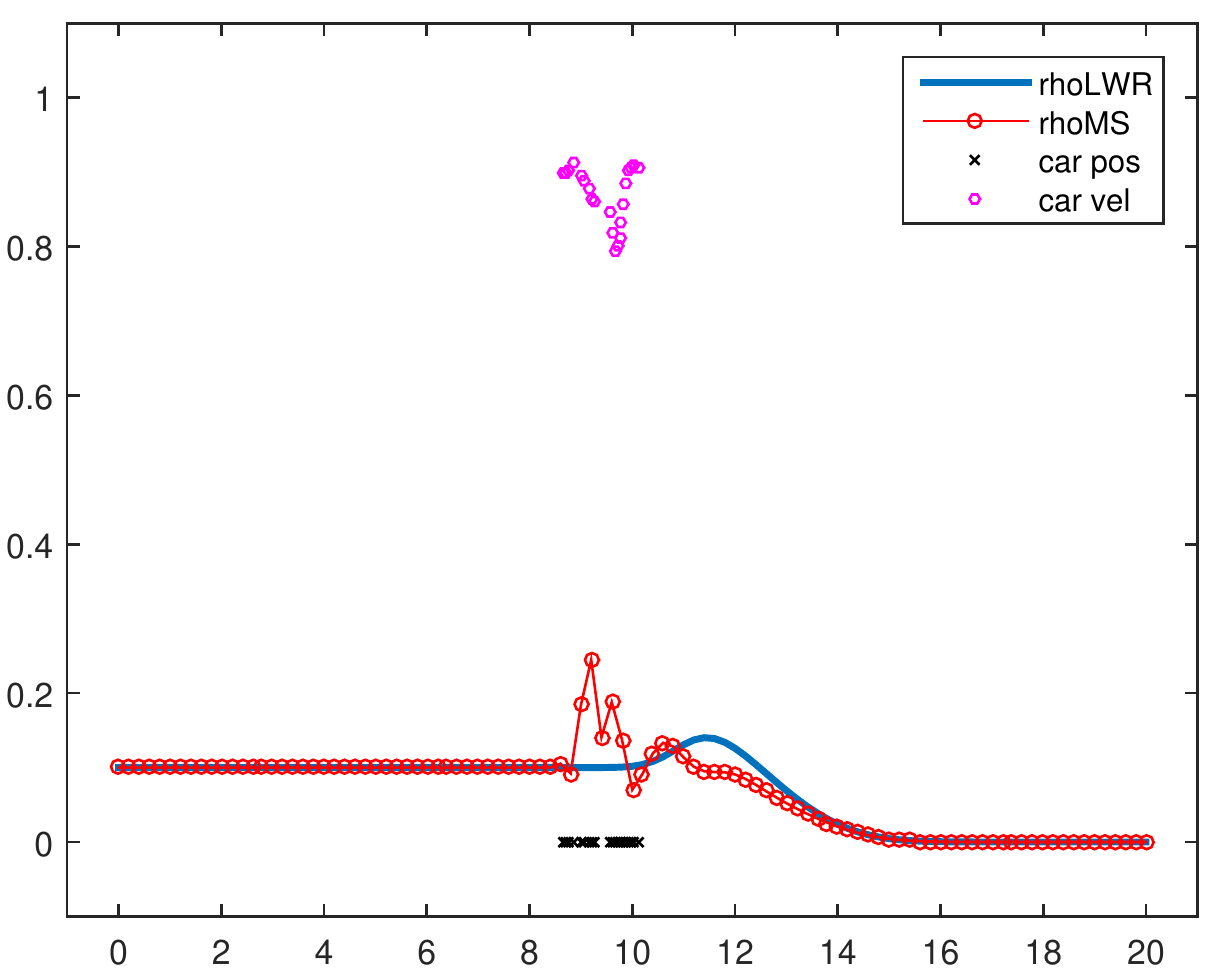}
		\put(52,-4){$x$}\put(-5,45){$\rho$}\end{overpic}
	\textbf{d.}
	\begin{overpic}[width=0.45\textwidth]{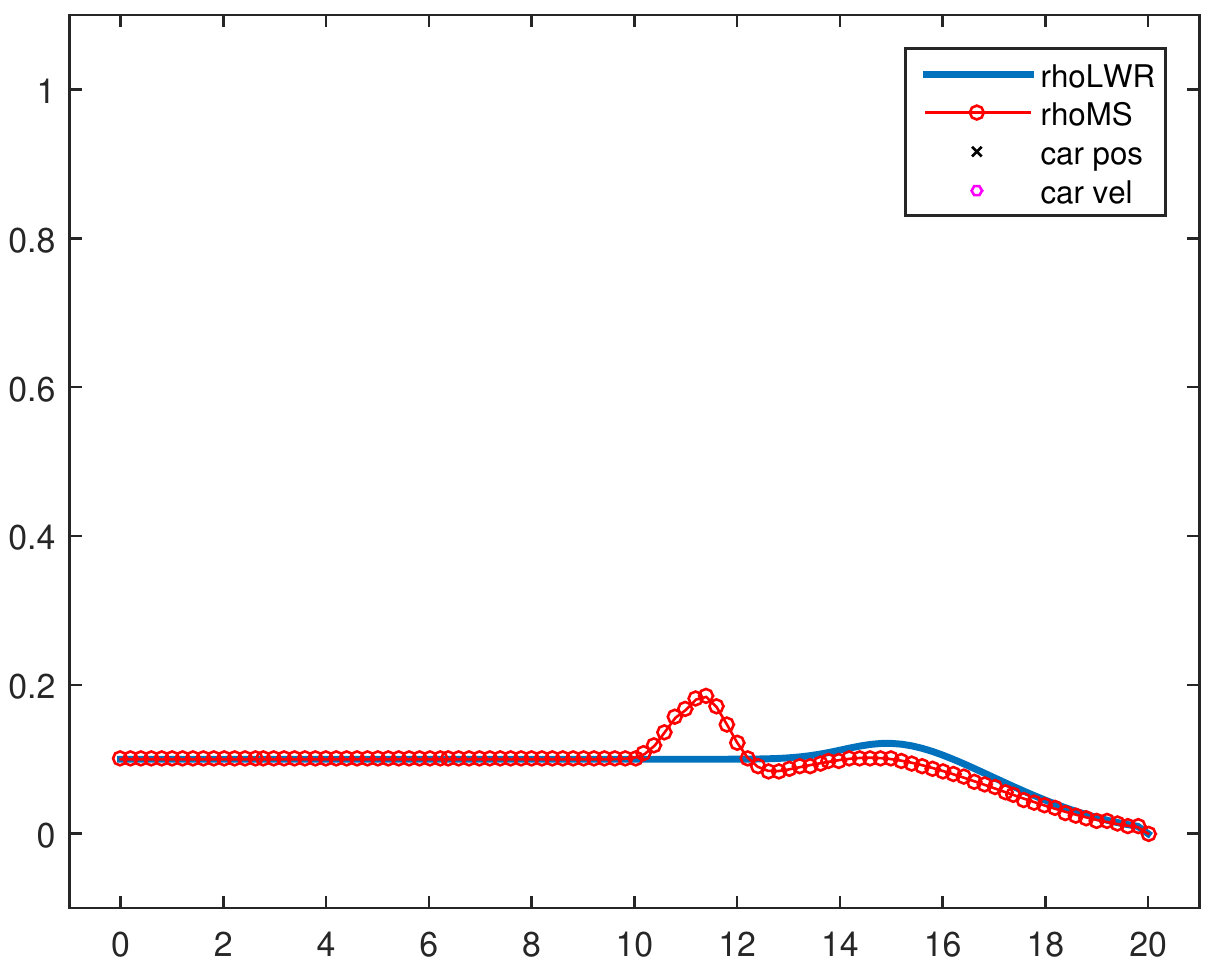}
		\put(52,-4){$x$}\put(-5,45){$\rho$}\end{overpic}
\caption{Test 3: \textbf{a.} $n=1$, \textbf{b.} $n=22$, \textbf{c.} $n=441$, \textbf{d.} $n=926$.}
\label{fig:T3}
\end{figure}
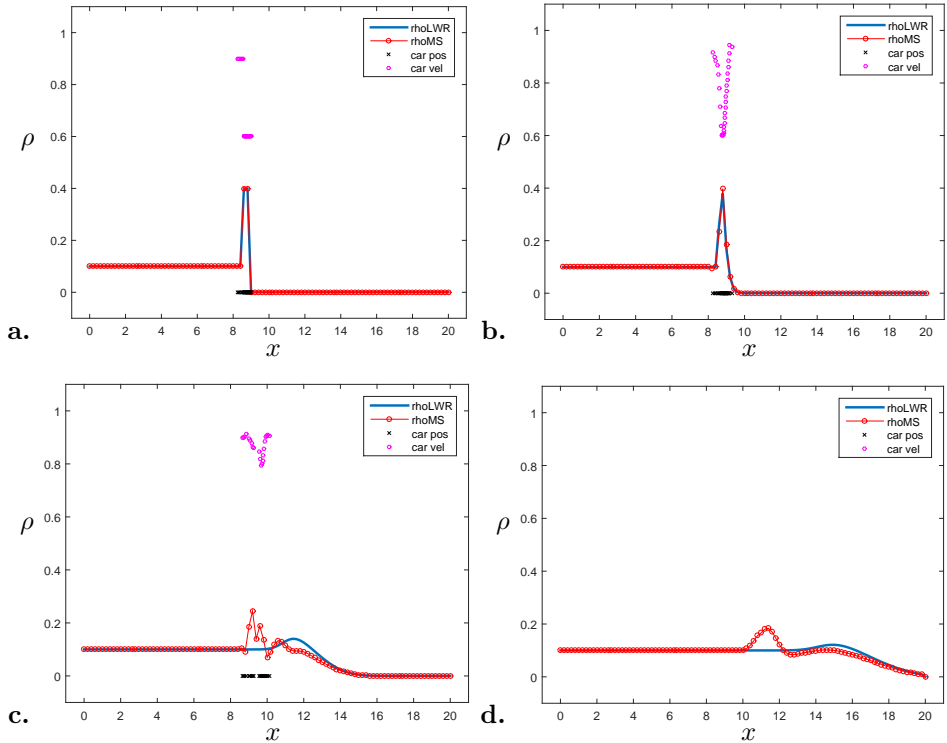

Microscopic vehicles are immediately activated in the region of perturbation (Fig.\ \ref{fig:T3}a), and, for a short time the multi-scale model and the LWR model behave similarly (Fig.\ \ref{fig:T3}b). At the microscopic level, it is clear that vehicles decelerates when get closer to the perturbation and then accelerate again (Fig.\ \ref{fig:T3}b).
After that, the LWR model tends to smear out the perturbation as usual. The multi-scale model, instead, \emph{self-sustains} the perturbation, which does not disappear, at least for a certain time (Fig.\ \ref{fig:T3}c). 
Around time step $n=900$ microscopic vehicles disappear since the perturbation is no longer strong enough to destroy the equilibrium, and the multi-scale model turns to be the LWR one (Fig.\ \ref{fig:T3}d).
%\begin{figure}[h!]
%\centering
%\textbf{a.}\includegraphics[width=0.47\textwidth]{./fig/T3_1.pdf} 
%\textbf{b.}\includegraphics[width=0.47\textwidth]{./fig/T3_22.pdf}\\
%\textbf{c.}\includegraphics[width=0.47\textwidth]{./fig/T3_441.pdf}
%\textbf{d.}\includegraphics[width=0.47\textwidth]{./fig/T3_926.pdf}

\subsection{Test 4 (stop \& go waves)}
\begin{figure}[h!]
	\centering
	\textbf{a.}
	\begin{overpic}[width=0.45\textwidth]{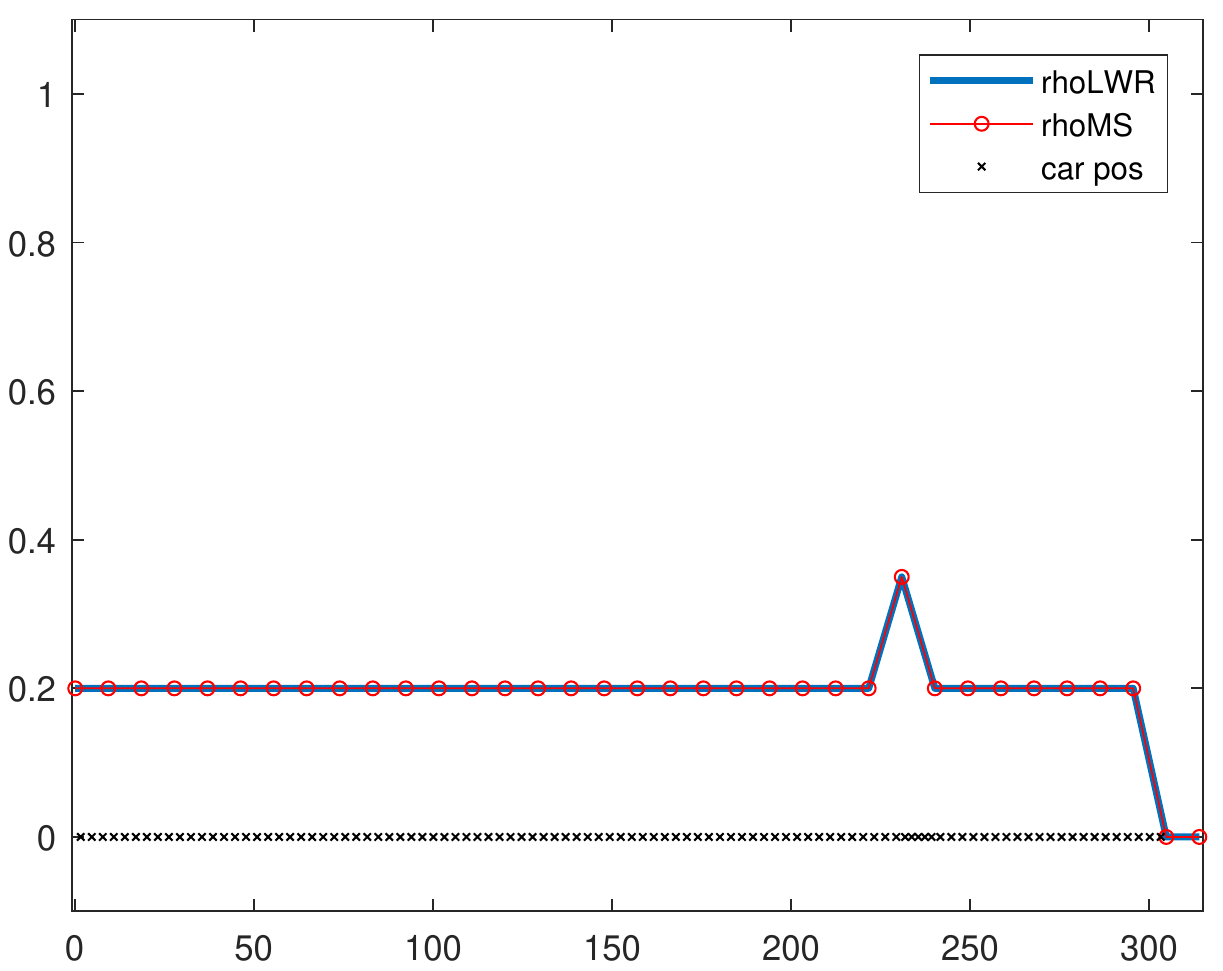}
		\put(52,-4){$x$}\put(-5,45){$\rho$}\end{overpic} 
	\textbf{b.}
	\begin{overpic}[width=0.45\textwidth]{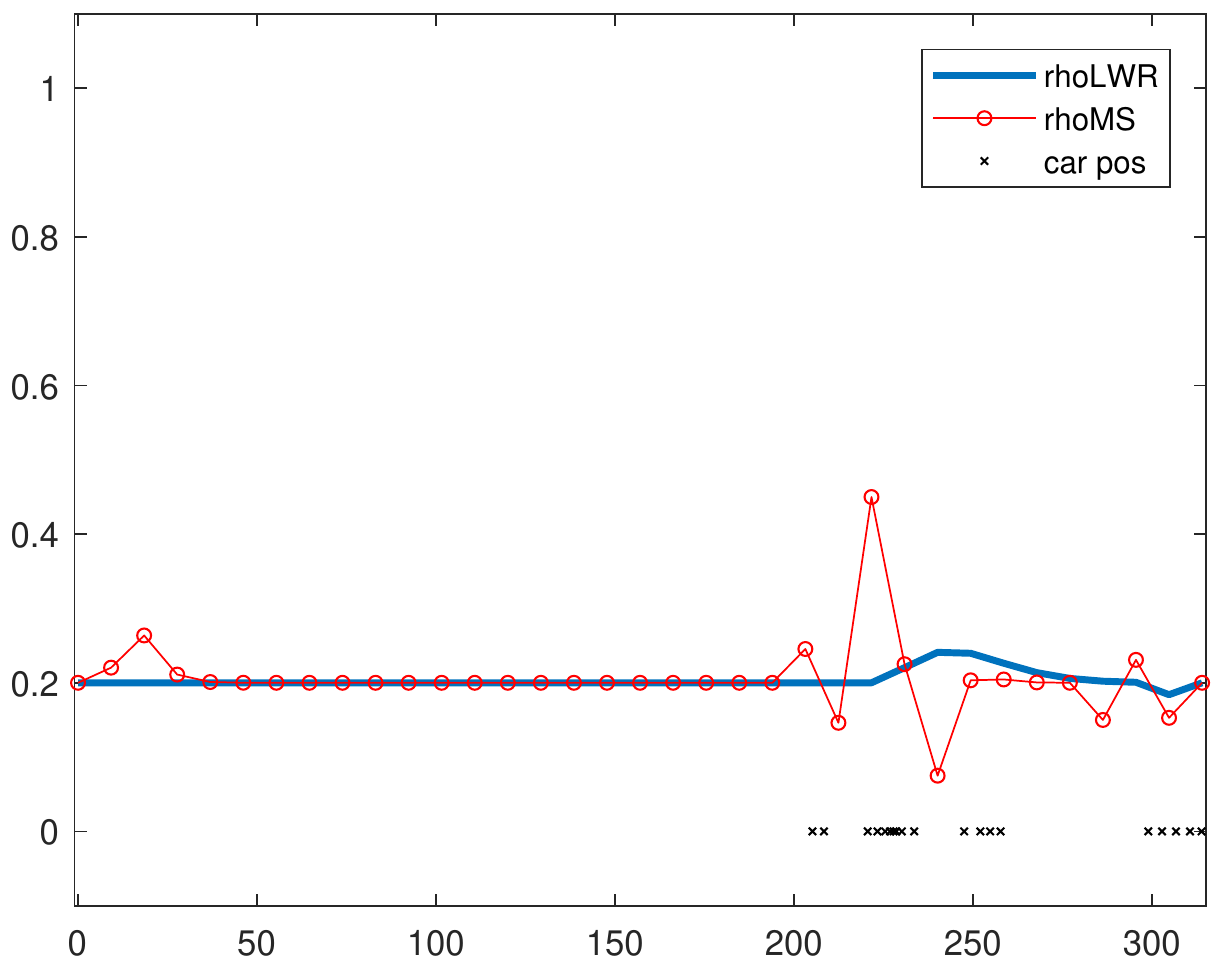}
		\put(52,-4){$x$}\put(-5,45){$\rho$}\end{overpic}
	\vskip0.5cm
	\textbf{c.}
	\begin{overpic}[width=0.45\textwidth]{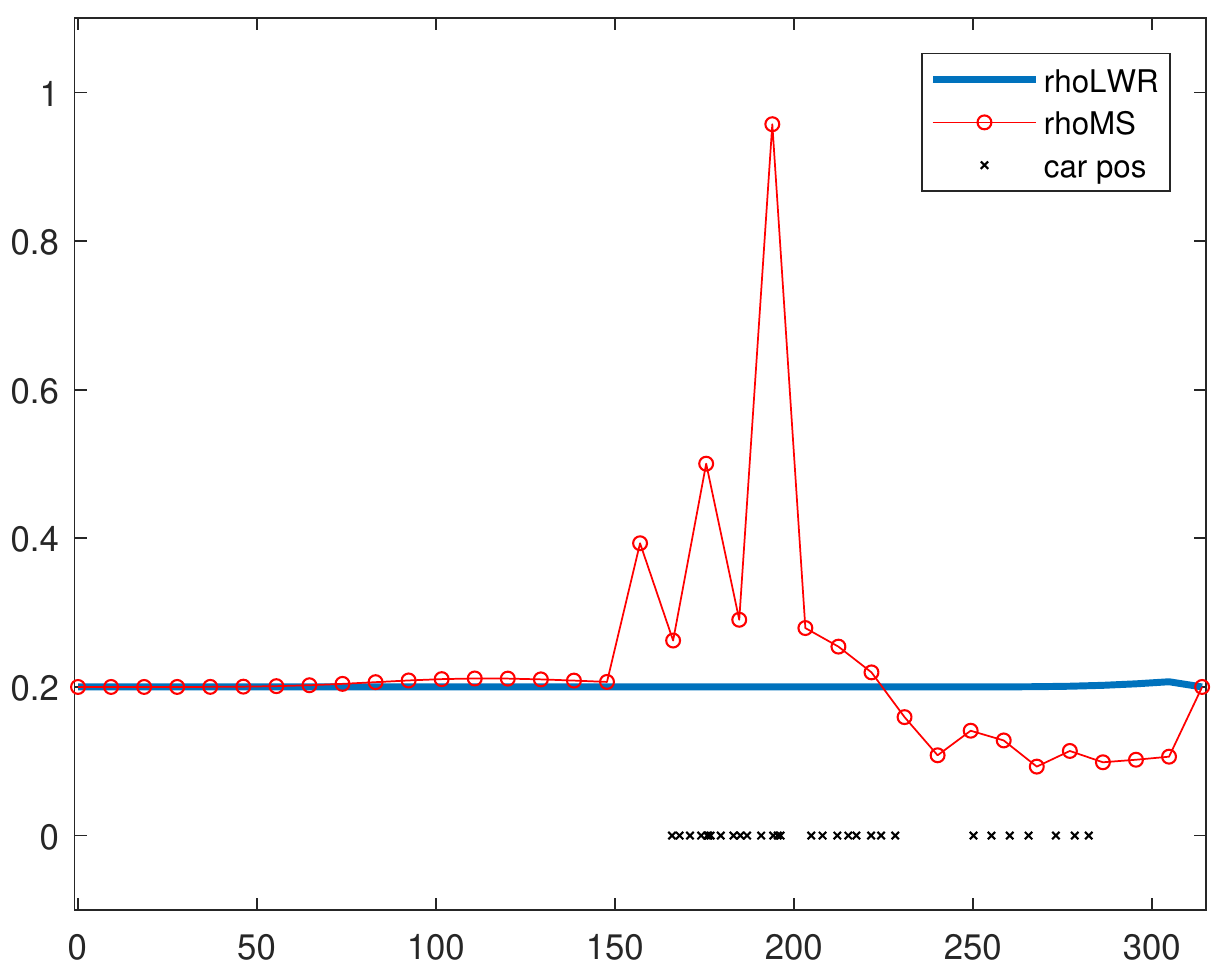}
		\put(52,-4){$x$}\put(-5,45){$\rho$}\end{overpic}
	\textbf{d.}
	\begin{overpic}[width=0.45\textwidth]{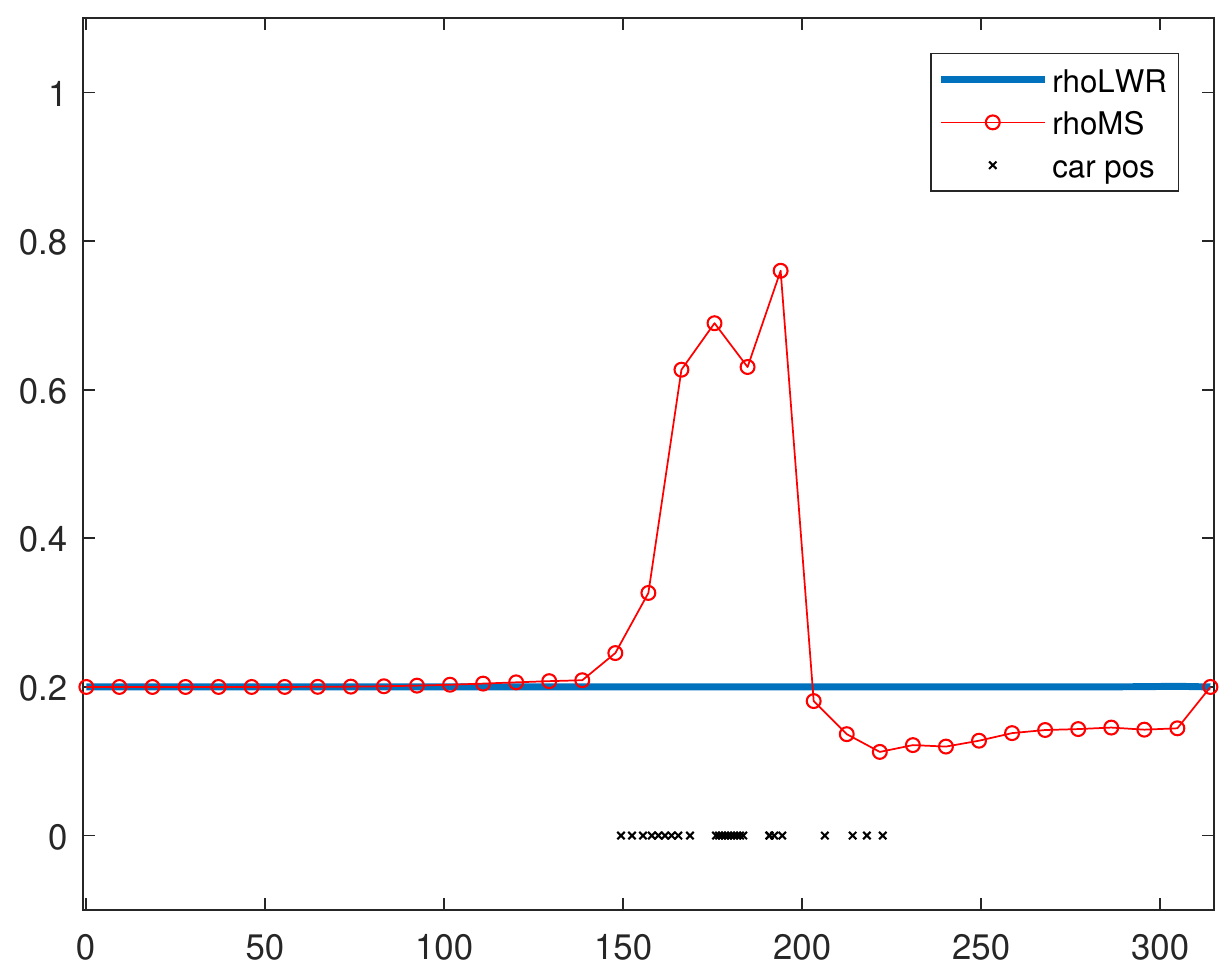}
		\put(52,-4){$x$}\put(-5,45){$\rho$}\end{overpic}
\caption{Test 4: \textbf{a.} $n=1$, \textbf{b.} $n=277$, \textbf{c.} $n=1666$, \textbf{d.} $n=2191$.}
\label{fig:T4}
\end{figure}
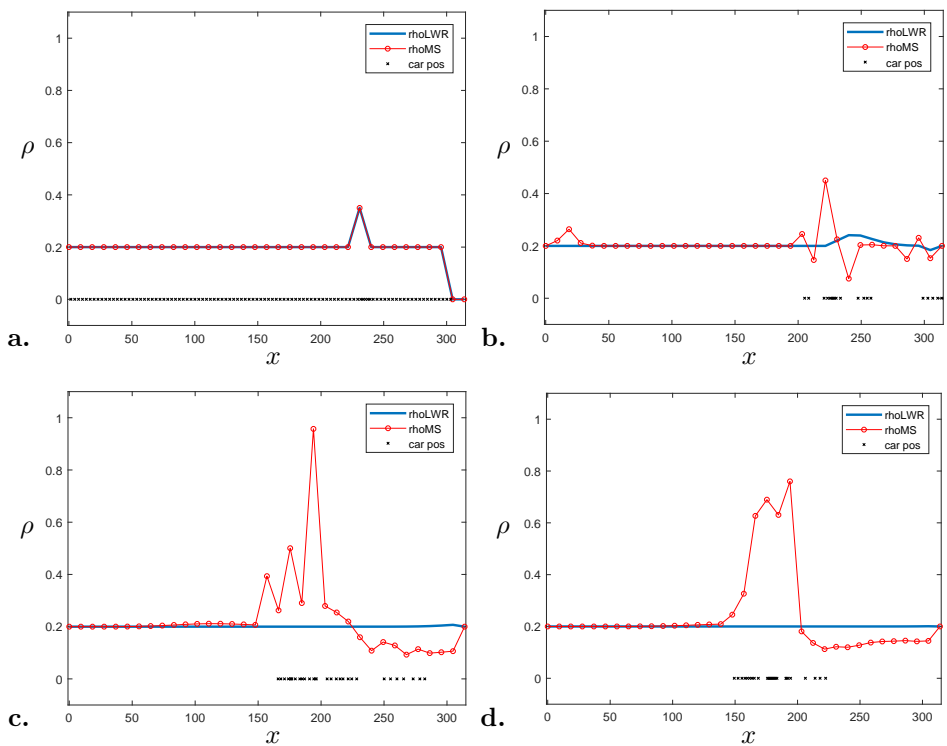

In this test we employ the model introduced in Section \ref{sec:cinese} specifically designed to reproduce stop \& go waves. Similarly to the previous test, a perturbation is created at initial time, but this time microscopic vehicles are forced to be activated everywhere in the domain (Fig.\ \ref{fig:T4}a). 
After time $\delta t$, vehicles located in regions at equilibrium are correctly deactivated, while they stay alive around the initial perturbation (Fig.\ \ref{fig:T4}b). 
After that, perturbation increases until vehicles almost stop completely (Fig.\ \ref{fig:T4}c), and finally a large stop \& go wave is formed (Fig.\ \ref{fig:T4}d).
%\begin{figure}[h!]
%\centering
%\textbf{a.}\includegraphics[width=0.47\textwidth]{./fig/T4_1.pdf} 
%\textbf{b.}\includegraphics[width=0.47\textwidth]{./fig/T4_277.pdf}\\
%\textbf{c.}\includegraphics[width=0.47\textwidth]{./fig/T4_1666.pdf}
%\textbf{d.}\includegraphics[width=0.47\textwidth]{./fig/T4_2191.pdf}

%
%
%
%
%
%

%\begin{figure}[htb]
%\centerline{
%\includegraphics[width=0.75\textwidth]{./fig/dftestcodatau01sr02.pdf} \,
%\includegraphics[width=0.75\textwidth]{./fig/dftestcodatau2sr02.pdf}
%}
%\caption{Fundamental diagram $\tau=0.1$ (left), $\tau=2$ (right).}
%\label{fig:dftestcoda}
%\end{figure}

\section*{Conclusion and future work}
In this paper we have proposed a new way to couple a microscopic and a macroscopic description in the framework of traffic flow modeling.
Exploiting this coupling, we have shown how to confine an expensive second-order (acceleration-based) microscopic description in space-time regions where bounded acceleration plays an important role in the dynamics. Regions where the traffic is moving at equilibrium can be instead safely described by simple and CPU-saving first-order macroscopic models.
The main limitation of the proposed approach is the presence of three additional parameters to be tuned, which rule the coupling.
Moreover, it remains to show the invariance of the interval $[0,\rhomax]$ in the numerical approximation, i.e.\ the fact that $\rho^0_j\in[0,\rhomax] \ \forall j \Rightarrow \rho_j^n\in[0,\rhomax] \ \forall j,n$. 
We claim that this property holds true only under suitable conditions on the discretization and model parameters.

Natural extensions of the proposed approach are given by nonconstant modeling parameters. For example, allowing $\theta$ to be dependent on $\rho$ itself, one can obtain a smooth passage from the macroscopic to the microscopic description as the traffic conditions become more and more congested. 
Also the relaxation parameter $\tau$ can depend on $\rho$, in such a way that drivers are described as more reactive if traffic conditions are more dangerous (e.g., in case of large flux characterized by a large number of vehicle moving at high velocity).

Another research direction could be the analytical study of the model \eqref{FtL2o}-\eqref{def:Acinese}-\eqref{def:v^ZZ}. Despite its simplicity, the model shows very interesting properties and features. One can understand under which conditions stop \& go waves arise and disappear and the role of the stochasticity in the original model. One can also investigate the many-particle limit of the model, together with the definition of the macroscopic counterpart of the parameter $\Delta_{\textup{min}}$.

\section*{Acknowledgments}
Authors want to thank Marco Di Francesco, Corrado Lattanzio and Andrea Tosin for the useful comments on the manuscript.

%\bibliographystyle{aims.bst}
%\bibliography{biblio_traffico}

%\medskip
% The data information below will be filled by AIMS editorial staff
%Received xxxx 20xx; revised xxxx 20xx.
%\medskip

\end{document}